\documentclass[10pt,oneside,reqno]{amsart}
\usepackage{amsfonts}
\usepackage{eurosym}
\usepackage{latexsym}
\usepackage{amssymb,amsmath,amsfont}
\usepackage[colorlinks=true]{hyperref}
\usepackage{graphicx}

\newtheorem{theorem}{Theorem}[section]
\newtheorem{lemma}[theorem]{Lemma}

\theoremstyle{definition}
\newtheorem{definition}[theorem]{Definition}

\newtheorem{example}[theorem]{Example}
\theoremstyle{remark}

\newcommand{\argmin}{\operatornamewithlimits{argmin}}
\numberwithin{equation}{section}
\global\parskip 6pt \oddsidemargin
0.1cm \evensidemargin 0.4cm \headheight 0.1cm \topmargin -2.0cm
\begin{document}
\title{Image restoration using an inertial viscosity fixed point algorithm}
\author{ $^{1}$ Ebru ALTIPARMAK AND $^{2}$ Ibrahim KARAHAN }
\address{$^{1}$ Department of Mathematics, Faculty of Science, Erzurum
Technical University, Erzurum, 25700, Turkey,}
\address{$^{2}$ Department of Mathematics, Faculty of Science, Erzurum
Technical University, Erzurum, 25700, Turkey,}
\email{$^{1}$ebru.altiparmak@erzurum.edu.tr,}
\email{$^{2}$ibrahimkarahan@erzurum.edu.tr,}

\begin{abstract}
The image restoration problem is one of the popular topics in image
processing studied by many authors on account of its applications in various
areas. The aim of this paper is to present a new algorithm by using
viscosity approximation with inertial effect for finding a common fixed
point of an infinite family of nonexpansive mappings in a Hilbert space and
obtaining more quality images from degenerate images. Some strong
convergence theorems are proved under mild conditions. The obtained results
are applied to solve monotone inclusion problems, convex minimization
problems, variational inequality problems and generalized equilibrium
problems. It is shown that the proposed algorithm performs better than some
other algorithms. Also, the effects of inertial and viscosity terms in the
algorithm on image restoration have been investigated.
\end{abstract}

\maketitle

\section{Introduction}

Image processing is a method to facilitate the perception of images by
computers and interpretation by humans. Images may be distorted for some
reason such as camera or object movement, electricity, heat, sharp and
sudden disturbances in the image signal and sensor illumination levels.
Image processing methods perform some operations on an image to obtain an
enhanced image or to extract some useful information from it. In recent
years, these methods have been used in almost every field such as military
industry, Forensic sciences, underwater imaging, astronomy, physics, art,
biomedical remote sensing applications, geographic sciences, image and data
storage, medical imaging, astronomical imaging and industrial automation.
One of the areas that have an important place in image processing is image
restoration. Image restoration is the process of obtaining a relatively
clear image from the distorted or noisy image. So, the goal of image
restoration techniques is to increase the quality of the images. After
reconstruction, the quality of the images can be measured with the values of
signal to noise ratio (SNR), improvement in signal to noise ratio (ISNR) and
peak signal to noise ratio (PSNR).The mathematical model for the image
restoration problem is formulated by%
\begin{equation*}
\upsilon =Ax+b
\end{equation*}%
where $x$ is the original image, $A$ is the blurring matrix, $b$ is the
additive noise and $v$ is the observed image. The aim of the image
restoration problem is to minimize additive noise $b$ by using the observed
image $v$. The main approach for this problem is to solve the regularized
least squares minimization problem given by:%
\begin{equation}
x^{\ast }=\argmin_{x}\left\{ \frac{1}{2}\left\Vert Ax-v\right\Vert
_{2}^{2}+\lambda K(x)\right\} ,  \label{99}
\end{equation}%
where $\lambda >0$ is a regularization parameter, $K(x)$ is a regularizer
function that should be convex and $\left\Vert Ax-v\right\Vert _{2}^{2}$ is
a least squares term that measures the distance between $h$ and $Ax$. In
1977, Tikhonov and Arsenin \cite{1} introduced the following Tikhonov
regularization problem by taking Tikhonov matrix $L$ as a special case of $K$%
:
\begin{equation*}
x^{\ast }=\argmin_{x}\left\{ \frac{1}{2}\left\Vert Ax-v\right\Vert
_{2}^{2}+\lambda \left\Vert Lx\right\Vert _{2}^{2}\right\} .
\end{equation*}%
By taking the $W$ wavelet transform matrix instead of the function $K$ in
problem (\ref{99}), we get the wavelet-based regularization given by:%
\begin{equation*}
x^{\ast }=\argmin_{x}\left\{ \frac{1}{2}\left\Vert Ax-v\right\Vert
_{2}^{2}+\lambda \left\Vert Wx\right\Vert _{1}\right\} .
\end{equation*}%
On the other hand, another successful regularization problem is known as $%
l_{1}$ regularization which is given as follows:
\begin{equation*}
x^{\ast }=\argmin_{x}\left\{ \frac{1}{2}\left\Vert Ax-v\right\Vert
_{2}^{2}+\lambda \left\Vert x\right\Vert _{1}\right\} .
\end{equation*}%
All these problems can be reformulated by in the following general way:%
\begin{equation}
x^{\ast }=\argmin_{x}\left\{ h\left( x\right) +g\left( x\right)
\right\} .  \label{102}
\end{equation}%
So, it is clear that the image restoration problem is a classical
minimization problem for the sum of two special functions. It is known that
the problem (\ref{102}) is equivalent to the following variational inclusion
problem:
\begin{equation}
0\in \nabla h\left( x^{\ast }\right) +\partial g\left( x^{\ast }\right) ,
\label{102a}
\end{equation}%
where $\partial g$ is the subdifferential of function $g$ defined by $%
\partial g\left( x\right) :=\left\{ u\in H:g\left( x\right) +\left\langle
y-x,u\right\rangle \leq g\left( y\right) ,\forall y\in H\right\} $ and $%
\nabla h$ is the gradient operator of function $h$\ and to the following
fixed point problem:%
\begin{equation}
x^{\ast }=J_{\lambda }^{\partial g}(I-\lambda \nabla f)x^{\ast },
\label{105}
\end{equation}%
where $J_{\lambda }^{\partial g}$ is the resolvent operator defined by $%
J_{\lambda }^{\partial g}=\left( I+\lambda \partial g\right) ^{-1}.$ Hence,
fixed point theory plays a very important role in solving image restoration
problems and so fixed point algorithms can be used to find the solutions of
the considered problems, see, for instance \cite{in, by, cho, ch, com, dong,
han, kun, kittt}. Lions and Mercier \cite{lion} introduced the following
classical forward backward splitting $\left( FBS\right) $ algorithm which is
one of the most important algorithm:
\begin{equation}
x_{n+1}=J_{\lambda }^{\partial g}(I-\lambda \nabla f)x_{n},  \label{60}
\end{equation}%
where $\lambda >0$ and $I$ is the identity operator. Lorenz and Pock \cite%
{lo} presented the inertial forward backward algorithm $\left( IFBS\right) $
for monotone operators in the following form:%
\begin{equation}
\left\{
\begin{array}{ll}
y_{n}=x_{n}+\theta _{n}\left( x_{n}-x_{n-1}\right) &  \\
x_{n+1}=J_{\lambda }^{\partial g}(I-\lambda \nabla f)y_{n}, & \forall n\geq 1%
\end{array}%
\right.  \label{61}
\end{equation}%
where $\theta _{n}$ is the inertial parameter which controls the momentum $%
x_{n}-x_{n-1}$. In 2016, Shehu and Cai \cite{choo} presented the following
algorithm by combining the algorithm (\ref{60}) with viscosity approximation
method:%
\begin{equation}
x_{n+1}=%
\begin{array}{cc}
\alpha _{n}f\left( x_{n}\right) +\left( 1-\alpha _{n}\right) J_{\lambda
_{n}}^{B}(I-\lambda _{n}A)x_{n}, & \forall n\geq 1%
\end{array}
\label{62}
\end{equation}%
where $f$ is a contraction mapping, $A$ is a $\nu $-inverse strongly
accretive mapping and $B$ is an $m$-accretive mapping. They proved the
convergence of the sequence generated by the algorithm in uniformly smooth
Banach spaces. In 2019, Kithuan et al. \cite{kit} proved some strong
convergence theorems for the following inertial viscosity forward-backward
splitting algorithm in a Hilbert space $H:$

\begin{equation}
\left\{
\begin{array}{l}
y_{n}=x_{n}+\theta _{n}\left( x_{n}-x_{n-1}\right) \\
x_{n+1}=\delta _{n}\nabla h\left( x_{n}\right) +\left( 1-\delta _{n}\right)
J_{\lambda _{n}}^{B}(I-\lambda _{n}A)y_{n},\forall n\geq 1%
\end{array}%
\right.  \label{63}
\end{equation}%
where $A:H\rightarrow H$ is a $\nu $-inverse strongly monotone mapping and $%
B:H\rightarrow 2^{H}$ is a maximal monotone operator. In 2017, Verma and
Shukla \cite{inertial} studied the new accelerated proximal gradient
algorithm $\left( NAGA\right) $ for a nonexpansive mapping in the following
way :
\begin{equation}
\left\{
\begin{array}{l}
y_{n}=x_{n}+\theta _{n}\left( x_{n}-x_{n-1}\right) \\
x_{n+1}=J_{\lambda _{n}}^{B}(I-\lambda _{n}A)\left( \left( 1-\beta
_{n}\right) y_{n}+\beta _{n}J_{\lambda _{n}}^{B}(I-\lambda
_{n}A)y_{n}\right) ,\forall n\geq 1.%
\end{array}%
\right.  \label{64}
\end{equation}

Padcharoen and Kumam \cite{kumam} introduced a modified MM algorithm $\left(
FBMMMA\right) $ for finding a common fixed point of a countable family of
nonexpansive operators in the following manner:%
\begin{equation}
\left\{
\begin{array}{l}
y_{n}=x_{n}+\theta _{n}\left( x_{n}-x_{n-1}\right) \\
z_{n}=\left( 1-\rho _{n}\right) y_{n}+\rho _{n}J_{\lambda
_{n}}^{B}(I-\lambda _{n}A)y_{n} \\
w_{n}=\left( 1-\delta _{n}-\rho _{n}\right) z_{n}+\delta _{n}J_{\lambda
_{n}}^{B}(I-\lambda _{n}A)z_{n}+\phi _{n}J_{\lambda _{n}}^{B}(I-\lambda
_{n}A)y_{n} \\
x_{n+1}=J_{\lambda _{n}}^{B}(I-\lambda _{n}A)w_{n}.%
\end{array}%
\right.  \label{65}
\end{equation}

They obtained the weak convergence of the sequence generated by the
algorithm in Hilbert space. Puangpee and Suantai \cite{sut} presented a new
accelerated fixed point algorithm $\left( AVFB\right) $ as follows :%
\begin{equation}
\left\{
\begin{array}{l}
y_{n}=x_{n}+\theta _{n}\left( x_{n}-x_{n-1}\right) \\
z_{n}=\left( 1-\sigma _{n}\right) y_{n}+\sigma _{n}J_{\lambda
_{n}}^{B}(I-\lambda _{n}A)y_{n} \\
x_{n+1}=\psi _{n}f\left( x_{n}\right) +\varrho _{n}J_{\lambda
_{n}}^{B}(I-\lambda _{n}A)y_{n}+\varphi _{n}J_{\lambda _{n}}^{B}(I-\lambda
_{n}A)z_{n}.%
\end{array}%
\right.  \label{66}
\end{equation}

They showed the strong convergence of the proposed algorithm for an infinite
family of nonexpansive mappings in Hilbert space. All these authors gave
some applications for image restoration problems and intended to obtain more
quality images.

In this paper, motivated and inspired by the given algorithms and $K$-
iteration algorithm introduced by Hussain et al. \cite{hussain}, we
introduced a new algorithm which is more effective than the algorithms exist
in the literature. We proved the strong convergence of generated sequence
and gave some application of the proposed algorithm to the different
problems especially to the image restoration problem\ to get better quality
images. This research is organized as follows. In section 2, Preliminaries,
we give some definitions and lemmas which we need to prove the main result.
In section 3, Main results, we prove the strong convergence of the proposed
algorithm. In the next section, Applications, we apply our main result to
solving inclusion problems, variational inequality problems, generalized
equilibrium problems and image restoration problems.

\section{Preliminaries}

Throughout this paper, let $H$ be a real Hilbert space with the inner
product $\left\langle .,.\right\rangle $ and the induced norm $\left\Vert
.\right\Vert $. Also, let $C$ be a nonempty closed and convex subset of a
real Hilbert space $H$ and $T$ a mapping on $C$. Then, the set of all fixed
points of $T$ is denoted by $F\left( T\right) :=\left\{ x\in C:x=Tx\right\}
. $

The metric projection $P_{C}:H\rightarrow C$ is defined as follows: by given
$x\in H$ there exist a unique point in $C$ such that%
\begin{equation*}
\left\Vert x-P_{C}x\right\Vert \leq \inf_{y\in C}\left\Vert x-y\right\Vert .
\end{equation*}

It is well-known that $P_{C}$ is a nonexpansive mapping and it can be
characterized by%
\begin{equation}
\left\langle x-P_{C}x,y-P_{C}x\right\rangle \leq 0  \label{30}
\end{equation}%
for all $y\in C,$ see \cite{gobell} for more details.

\begin{definition}
\cite{bau} Let $T:C\rightarrow H$ \ be a nonlinear operator. Then $T$ is
said to be:

\begin{enumerate}
\item[(1)] $L$-Lipschitz continuous, if there exists a constant $L>0$ such
that:%
\begin{equation*}
\left\Vert Tx-Ty\right\Vert \leq L\left\Vert x-y\right\Vert ,\forall x,y\in
C.
\end{equation*}

If $L=1,$ then $T$ is said to be nonexpansive mapping and if $L<1$ then $T$
is said to be contraction mapping.

\item[(2)] $\upsilon $-strongly monotone operator if there exists $\upsilon
>0$ such that:%
\begin{equation*}
\left\langle Tx-Ty,x-y\right\rangle \geq \upsilon \left\Vert x-y\right\Vert
^{2},\text{ }\forall x,y\in C
\end{equation*}

\item[(3)] $\eta $-inverse strongly monotone operator if there exists $\eta
>0$ such that:%
\begin{equation*}
\left\langle Tx-Ty,x-y\right\rangle \geq \eta \left\Vert Tx-Ty\right\Vert
^{2},\text{ }\forall x,y\in C.
\end{equation*}
\end{enumerate}
\end{definition}

We can see that if $T$ is $\eta $-inverse strongly monotone, then it is $%
\frac{1}{\eta }$-Lipschitz continuous.

Let $A:H\rightarrow 2^{H}$ be a set-valued operator. $A$ is called\ monotone
if $\left\langle z-w,x-y\right\rangle \geq 0,$ for all $z\in Ax$ and $w\in
Ay.$ If the graph of a monotone operator is not properly contained in the
graph of any other monotone operators, then it is called the maximal
monotone operator. The following lemmas give some useful informations
regarding maximal monotonicity.

\begin{lemma}
\cite{bau}\label{minty} Let $A:H\rightarrow 2^{H}$ be a monotone operator.
Then $A$ is maximal monotone if and only if $R\left( I+A\right) =H.$
\end{lemma}

\begin{lemma}
\cite{bau} Let $\Gamma _{0}\left( H\right) $ denotes the class of all lower
semi-continuous convex functions defined from $H$ to $\left( -\infty ,\infty %
\right] .$ If $g\in \Gamma _{0}\left( H\right) ,$ then $\partial g$ is
maximal monotone.
\end{lemma}

Let $A:H\rightarrow 2^{H}$ be a maximal monotone operator. Then the
resolvent operator $J_{\lambda }^{A}:H\rightarrow H$ \ associated with $A$
is defined by $J_{\lambda }^{A}=\left( I+\lambda A\right) ^{-1}$ for all $%
x\in H$ and for some $\lambda >0.$ It is well-known that $J_{\lambda }^{A}$
is a nonexpansive mapping and $F\left( J_{\lambda }^{A}\right) =A^{-1}0$
where $A^{-1}0=\left\{ x\in H:0\in Ax\right\} $ which is called the set of
all zero points of $A.$

\begin{definition}
\cite{bau} Let $g\in $ $\Gamma _{0}\left( H\right) $ and $\lambda >0.$ The
proximal operator of parameter $\lambda $ of $g$ at $x$ is defined by%
\begin{equation*}
prox_{_{\lambda }g}\left( x\right) =\argmin_{y\in H}\left\{
g\left( y\right) +\frac{1}{2\lambda }\left\Vert y-x\right\Vert ^{2}\right\} .
\end{equation*}
\end{definition}

It is known that if $g\in $ $\Gamma _{0}\left( H\right) ,$ then $J_{\lambda
}^{\partial g}=prox_{\lambda g}.$ By taking $l_{1}$-norm instead of $g$ in
the definition of proximal operator, the proximal operator is called the
soft thresholding operator and it can be given by the following way:
\begin{equation*}
prox_{\lambda \left\Vert .\right\Vert _{1}}\left( x\right) =sign(x)\max
\left\{ \left\Vert x\right\Vert _{1}-\lambda ,0\right\} .
\end{equation*}

We recall the following properties which is useful to prove our main result,
please see for details \cite{nan, tak}.

Let $\left\{ T_{n}\right\} $ and $\Lambda $ be two families of nonexpansive
mappings of $H$ into itself such that $\emptyset \neq F\left( \Lambda
\right) \subset \cap _{n=1}^{\infty }F\left( T_{n}\right) $ where $F\left(
\Lambda \right) $ is the set of all common fixed points of mappings belongs
to $\Lambda .$Then we say that $\left\{ T_{n}\right\} $ satisfies
NST-condition (I) with $\Lambda $ if for every bounded sequence $\left\{
x_{n}\right\} $,
\begin{equation*}
\lim_{n\rightarrow \infty }\left\Vert x_{n}-T_{n}x_{n}\right\Vert =0\text{
implies }\lim_{n\rightarrow \infty }\left\Vert x_{n}-Tx_{n}\right\Vert =0,%
\text{ for all }T\in \Lambda .
\end{equation*}

Especially, if $\Lambda $ consist of one mapping, that is, $\Lambda =\left\{
T\right\} ,$ then $\left\{ T_{n}\right\} $ is said to satisfy the
NST-condition (I) with $T.$

\begin{lemma}
\label{FORWARD}\cite{bussaban} For a real Hilbert space $H$, let $%
A:H\rightarrow H$ be a $L$-Lipschitz operator and let $B:H\rightarrow 2^{H}$%
be a maximal monotone operator. If $T_{n}$ is a forward backward operator,
i.e. $T_{n}=J_{\lambda _{n}}^{B}(I-\lambda _{n}A)$ where $\lambda _{n}\in
\left( 0,2/L\right) $ for all $n\geq 1$ such that $\lambda _{n}$ converges
to $\lambda $, then $\left\{ T_{n}\right\} $ satisfies NST-condition (I)
with $T$ , where $T=J_{\lambda }^{B}(I-\lambda A)$ is a forward backward
operator with $\lambda \in \left( 0,2/L\right) .$
\end{lemma}

The following lemmas are important for proving our main theorem.

\begin{lemma}
\label{lemmaaa}\cite{gobel} Let $T$ be a nonexpansive operator from $H$ into
itself with $F\left( T\right) \neq \emptyset .$ Then, the mapping $I-T$ is
demiclosed at zero, that is, for any sequences $\left\{ x_{n}\right\} \in $ $%
H$ such that $x_{n}\rightharpoonup x$ $\in H$ and $\left\Vert
x_{n}-Tx_{n}\right\Vert \rightarrow 0$ as $n\rightarrow \infty ,$ then it
implies $x\in F\left( T\right) .$
\end{lemma}

\begin{lemma}
\label{a2}\cite{bau} Let $H$ be a real Hilbert space. Then for all \ $x,y\in
H,$ and $\lambda \in \left[ 0,1\right] ,$ the following properties hold:

\begin{enumerate}
\item[(1)] $\left\Vert x\pm y\right\Vert ^{2}=\left\Vert x\right\Vert
^{2}\pm 2\left\langle x,y\right\rangle +\left\Vert y\right\Vert ^{2},$

\item[(2)] $\left\Vert x+y\right\Vert ^{2}\leq \left\Vert x\right\Vert
^{2}+2\left\langle y,x+y\right\rangle ,$

\item[(3)] $\left\Vert \lambda x+\left( 1-\lambda \right) y\right\Vert
^{2}=\lambda \left\Vert x\right\Vert ^{2}+\left( 1-\lambda \right)
\left\Vert y\right\Vert ^{2}-\lambda \left( 1-\lambda \right) \left\Vert
x-y\right\Vert ^{2}.$
\end{enumerate}
\end{lemma}

\begin{lemma}
\label{lemmaa}\cite{aoya,xu} Let $\left\{ s_{n}\right\} $ and $\left\{
\varepsilon _{n}\right\} $ be sequences of nonnegative real numbers such
that
\begin{equation*}
s_{n+1}\leq \left( 1-\delta _{n}\right) s_{n}+\delta _{n}r_{n}+\varepsilon
_{n},
\end{equation*}

where $\left\{ \delta _{n}\right\} $ is a sequence in $\left[ 0,1\right] $
and $\left\{ r_{n}\right\} $ is a real sequence. If the following conditions
hold, then $\lim_{n\rightarrow \infty }s_{n}=0:$

\begin{enumerate}
\item[(1)] $\sum\limits_{n=1}^{\infty }\delta _{n}=\infty ,$

\item[(2)] $\sum\limits_{n=1}^{\infty }\varepsilon _{n}<\infty ,$

\item[(3)] $\limsup_{n\rightarrow \infty }r_{n}\leq 0.$
\end{enumerate}
\end{lemma}

\begin{lemma}
\label{main}\cite{main} Let $\left\{ \Phi _{n}\right\} $ be a sequence of
real numbers that does not decrease at infinity such that there exists a
subsequence $\left\{ \Phi _{n_{i}}\right\} $ of $\left\{ \Phi _{n}\right\} $
which satisfies $\Phi _{n_{i}}<\Phi _{n_{i+1}}$ for all $i\in
\mathbb{N}
.$ Let $\left\{ \tau \left( n\right) \right\} _{n\geq n_{0}}$ be a sequence
of integer, defined as follows:%
\begin{equation*}
\tau \left( n\right) :=\max \left\{ l\leq n:\Phi _{l}<\Phi _{l+1}\right\} .
\end{equation*}

Then the followings are satisfied:

\begin{enumerate}
\item[(1)] $\tau \left( n_{0}\right) \leq \tau \left( n_{0}+1\right) $ $\leq
...$ and $\tau \left( n\right) \rightarrow \infty ,$

\item[(2)] $\Phi _{\tau \left( n\right) }\leq \Phi _{\tau \left( n\right) +1%
\text{ }}$and $\Phi _{n}\leq \Phi _{\tau \left( n\right) +1\text{ }},$ for
all $n\geq n_{0}.$
\end{enumerate}
\end{lemma}

\section{Main Results}

In this section, we present a new algorithm for finding a common fixed point
of an infinite family of nonexpansive mappings in real Hilbert space and
prove its strong convergence under some mild assumptions.

\begin{theorem}
\label{theorem}\noindent\ Let $\left\{ T_{n}\right\} $ be a family of
nonexpansive mappings on $H$ into itself which satisfies the NST-condition
(I) with a nonexpansive mapping $T:H\rightarrow H$. Let $x_{0},x_{1}\in H$, $%
f:H\rightarrow H$ be a $k$-contraction mapping and $\left\{ x_{n}\right\} $
be a sequence generated by%
\begin{equation}
\left\{
\begin{array}{l}
z_{n}=x_{n}+\theta _{n}\left( x_{n}-x_{n-1}\right) \\
y_{n}=\left( 1-\beta _{n}\right) z_{n}+\beta _{n}T_{n}z_{n} \\
\omega _{n}=T_{n}\left( \left( 1-\alpha _{n}\right) T_{n}x_{n}+\alpha
_{n}T_{n}y_{n}\right) \\
x_{n+1}=\left( 1-\gamma _{n}\right) T_{n}\omega _{n}+\gamma _{n}f\left(
\omega _{n}\right)%
\end{array}%
\right.  \label{1''}
\end{equation}

where,%
\begin{equation*}
\theta _{n}:=\left\{
\begin{array}{ll}
\min \left\{ \theta ,\frac{\eta _{n}\gamma _{n}}{\left\Vert
x_{n}-x_{n-1}\right\Vert }\right\} & if\text{ \ }x_{n}\neq x_{n-1} \\
\theta , & otherwise%
\end{array}%
\right. ,\text{ }
\end{equation*}%
for $\theta \geq 0$, $\left\{ \eta _{n}\right\} \in \left( 0,\infty \right) $
and $\left\{ \alpha _{n}\right\} ,\left\{ \beta _{n}\right\} ,\left\{ \gamma
_{n}\right\} \in \left( 0,1\right) $ be sequences which satisfy the
conditions:

\begin{enumerate}
\item[(1)] $0<a\leq \alpha _{n}<a^{^{\prime }}<1,$

\item[(2)] $\lim_{n\rightarrow \infty }\eta _{n}=0,$

\item[(3)] $\lim_{n\rightarrow \infty }\gamma _{n}=0,$ $\sum\limits_{n=1}^{%
\infty }\gamma _{n}=\infty ,$
\end{enumerate}

for some positive real numbers $a$ and $a^{^{\prime }}$. Then the sequence $%
\left\{ x_{n}\right\} $ converges strongly to a point $x^{\ast }$of $F\left(
T\right) ,$ where $x^{\ast }=P_{F\left( T\right) }f\left( x^{\ast }\right) .$
\end{theorem}

\begin{proof}
First, we prove that $\left\{ x_{n}\right\} $ is bounded. Let $x^{\ast }\in
F\left( T\right) $ such that $x^{\ast }=P_{F\left( T\right) }f\left( x^{\ast
}\right) $. By Algorithm \ref{1''}, we can write \
\begin{eqnarray}
\left\Vert z_{n}-x^{\ast }\right\Vert &=&\left\Vert x_{n}+\theta _{n}\left(
x_{n}-x_{n-1}\right) -x^{\ast }\right\Vert  \notag \\
&\leq &\left\Vert x_{n}-x^{\ast }\right\Vert +\theta _{n}\left\Vert
x_{n}-x_{n-1}\right\Vert ,  \label{1}
\end{eqnarray}

and, since $T_{n}$ is a nonexpansive mapping, we have
\begin{eqnarray}
\left\Vert y_{n}-x^{\ast }\right\Vert &=&\left\Vert \left( 1-\beta
_{n}\right) z_{n}+\beta _{n}T_{n}z_{n}-x^{\ast }\right\Vert  \notag \\
&\leq &\left( 1-\beta _{n}\right) \left\Vert z_{n}-x^{\ast }\right\Vert
+\beta _{n}\left\Vert T_{n}z_{n}-x^{\ast }\right\Vert  \notag \\
&=&\left( 1-\beta _{n}\right) \left\Vert z_{n}-x^{\ast }\right\Vert +\beta
_{n}\left\Vert T_{n}z_{n}-T_{n}x^{\ast }\right\Vert  \notag \\
&\leq &\left\Vert z_{n}-x^{\ast }\right\Vert ,  \label{2}
\end{eqnarray}

and also,
\begin{eqnarray}
\left\Vert \omega _{n}-x^{\ast }\right\Vert &=&\left\Vert T_{n}\left( \left(
1-\alpha _{n}\right) T_{n}x_{n}+\alpha _{n}T_{n}y_{n}\right) -x^{\ast
}\right\Vert  \notag \\
&=&\left\Vert T_{n}\left( \left( 1-\alpha _{n}\right) T_{n}x_{n}+\alpha
_{n}T_{n}y_{n}\right) -T_{n}x^{\ast }\right\Vert  \notag \\
&\leq &\left\Vert \left( 1-\alpha _{n}\right) T_{n}x_{n}+\alpha
_{n}T_{n}y_{n}-x^{\ast }\right\Vert  \notag \\
&=&\left\Vert \left( 1-\alpha _{n}\right) \left( T_{n}x_{n}-x^{\ast }\right)
+\alpha _{n}\left( T_{n}y_{n}-x^{\ast }\right) \right\Vert  \notag \\
&=&\left\Vert \left( 1-\alpha _{n}\right) \left( T_{n}x_{n}-T_{n}x^{\ast
}\right) +\alpha _{n}\left( T_{n}y_{n}-T_{n}x^{\ast }\right) \right\Vert
\notag \\
&\leq &\left( 1-\alpha _{n}\right) \left\Vert x_{n}-x^{\ast }\right\Vert
+\alpha _{n}\left\Vert y_{n}-x^{\ast }\right\Vert .  \label{3}
\end{eqnarray}

Combining (\ref{1}), (\ref{2}) and (\ref{3}), we obtain that%
\begin{eqnarray}
\left\Vert x_{n+1}-x^{\ast }\right\Vert &=&\left\Vert \left( 1-\gamma
_{n}\right) T_{n}\omega _{n}+\gamma _{n}f\left( \omega _{n}\right) -x^{\ast
}\right\Vert  \notag \\
&=&\left\Vert \left( 1-\gamma _{n}\right) \left( T_{n}\omega _{n}-x^{\ast
}\right) +\gamma _{n}\left( f\left( \omega _{n}\right) -f\left( x^{\ast
}\right) \right) +\gamma _{n}\left( f\left( x^{\ast }\right) -x^{\ast
}\right) \right\Vert  \notag \\
&\leq &\left( 1-\gamma _{n}\right) \left\Vert T_{n}\omega _{n}-x^{\ast
}\right\Vert +\gamma _{n}\left\Vert f\left( \omega _{n}\right) -f\left(
x^{\ast }\right) \right\Vert +\gamma _{n}\left\Vert f\left( x^{\ast }\right)
-x^{\ast }\right\Vert  \notag \\
&=&\left( 1-\gamma _{n}\right) \left\Vert T_{n}\omega _{n}-T_{n}x^{\ast
}\right\Vert +\gamma _{n}\left\Vert f\left( \omega _{n}\right) -f\left(
x^{\ast }\right) \right\Vert +\gamma _{n}\left\Vert f\left( x^{\ast }\right)
-x^{\ast }\right\Vert  \notag \\
&\leq &\left( 1-\gamma _{n}\right) \left\Vert \omega _{n}-x^{\ast
}\right\Vert +\gamma _{n}k\left\Vert \omega _{n}-x^{\ast }\right\Vert
+\gamma _{n}\left\Vert f\left( x^{\ast }\right) -x^{\ast }\right\Vert  \notag
\\
&\leq &\left( 1-\gamma _{n}\left( 1-k\right) \right) \left[ \left( 1-\alpha
_{n}\right) \left\Vert x_{n}-x^{\ast }\right\Vert +\alpha _{n}\left\Vert
y_{n}-x^{\ast }\right\Vert \right] +\gamma _{n}\left\Vert f\left( x^{\ast
}\right) -x^{\ast }\right\Vert  \notag \\
&\leq &\left( 1-\gamma _{n}\left( 1-k\right) \right) \left[ \left( 1-\alpha
_{n}\right) \left\Vert x_{n}-x^{\ast }\right\Vert +\alpha _{n}\left\Vert
x_{n}-x^{\ast }\right\Vert \right.  \notag \\
&&\left. +\alpha _{n}\theta _{n}\left\Vert x_{n}-x_{n-1}\right\Vert \right]
+\gamma _{n}\left\Vert f\left( x^{\ast }\right) -x^{\ast }\right\Vert  \notag
\\
&\leq &\left( 1-\gamma _{n}\left( 1-k\right) \right) \left\Vert
x_{n}-x^{\ast }\right\Vert +a^{^{\prime }}\gamma _{n}.\frac{\theta _{n}}{%
\gamma _{n}}\left\Vert x_{n}-x_{n-1}\right\Vert +\gamma _{n}\left\Vert
f\left( x^{\ast }\right) -x^{\ast }\right\Vert .  \label{4}
\end{eqnarray}

By using the definition of $\theta _{n}$ and condition $\left( 2\right) $,
it is clear that%
\begin{equation*}
\frac{\theta _{n}}{\gamma _{n}}\left\Vert x_{n}-x_{n-1}\right\Vert
\rightarrow 0\text{ as }n\rightarrow \infty .
\end{equation*}

Hence, there exists a positive constant $M_{1}>0$ such that, for all $n\geq
1 $,%
\begin{equation*}
\frac{\theta _{n}}{\gamma _{n}}\left\Vert x_{n}-x_{n-1}\right\Vert \leq
M_{1}.
\end{equation*}

It follows from (\ref{4}) that,
\begin{eqnarray*}
\left\Vert x_{n+1}-x^{\ast }\right\Vert &\leq &\left( 1-\gamma _{n}\left(
1-k\right) \right) \left\Vert x_{n}-x^{\ast }\right\Vert +\gamma _{n}\left(
a^{^{\prime }}M_{1}+\left\Vert f\left( x^{\ast }\right) -x^{\ast
}\right\Vert \right) \\
&=&\left( 1-\gamma _{n}\left( 1-k\right) \right) \left\Vert x_{n}-x^{\ast
}\right\Vert +\gamma _{n}\left( 1-k\right) \left[ \frac{a^{^{\prime
}}M_{1}+\left\Vert f\left( x^{\ast }\right) -x^{\ast }\right\Vert }{\left(
1-k\right) }\right] \\
&\leq &\max \left\{ \left\Vert x_{n}-x^{\ast }\right\Vert ,\frac{a^{^{\prime
}}M_{1}+\left\Vert f\left( x^{\ast }\right) -x^{\ast }\right\Vert }{\left(
1-k\right) }\right\} \\
&&\vdots \\
&\leq &\max \left\{ \left\Vert x_{1}-x^{\ast }\right\Vert ,\frac{a^{^{\prime
}}M_{1}+\left\Vert f\left( x^{\ast }\right) -x^{\ast }\right\Vert }{\left(
1-k\right) }\right\} ,
\end{eqnarray*}

for all $n\geq 1.$ So, we obtain that $\left\{ x_{n}\right\} $ is bounded
and hence $\left\{ z_{n}\right\} ,$ $\left\{ y_{n}\right\} $ and $\left\{
\omega _{n}\right\} $ are also bounded.

Secondly, we want to prove that $x_{n}\rightarrow x^{\ast }=P_{F\left(
T\right) }f\left( x^{\ast }\right) .$ Indeed, we have the followings for all%
\begin{eqnarray}
\left\Vert z_{n}-x^{\ast }\right\Vert ^{2} &=&\left\Vert x_{n}+\theta
_{n}\left( x_{n}-x_{n-1}\right) -x^{\ast }\right\Vert ^{2}  \notag \\
&\leq &\left\Vert x_{n}-x^{\ast }\right\Vert ^{2}+2\theta _{n}\left\Vert
x_{n}-x^{\ast }\right\Vert \left\Vert x_{n}-x_{n-1}\right\Vert +\theta
_{n}^{2}\left\Vert x_{n}-x_{n-1}\right\Vert ^{2},  \label{5}
\end{eqnarray}

and%
\begin{eqnarray*}
\left\Vert \omega _{n}-x^{\ast }\right\Vert ^{2} &=&\left\Vert T_{n}\left(
\left( 1-\alpha _{n}\right) T_{n}x_{n}+\alpha _{n}T_{n}y_{n}\right) -x^{\ast
}\right\Vert ^{2} \\
&\leq &\left\Vert \left( 1-\alpha _{n}\right) T_{n}x_{n}+\alpha
_{n}T_{n}y_{n}-x^{\ast }\right\Vert ^{2} \\
&=&\left\Vert \left( 1-\alpha _{n}\right) \left( T_{n}x_{n}-T_{n}x^{\ast
}\right) +\alpha _{n}\left( T_{n}y_{n}-T_{n}x^{\ast }\right) \right\Vert
^{2}.
\end{eqnarray*}%
From (\ref{2}), (\ref{5}) and the property $(3)$ of Lemma (\ref{a2}), we get%
\begin{eqnarray}
\left\Vert \omega _{n}-x^{\ast }\right\Vert ^{2} &\leq &\alpha
_{n}\left\Vert T_{n}y_{n}-T_{n}x^{\ast }\right\Vert ^{2}+\left( 1-\alpha
_{n}\right) \left\Vert T_{n}x_{n}-T_{n}x^{\ast }\right\Vert ^{2}  \notag \\
&&-\alpha _{n}\left( 1-\alpha _{n}\right) \left\Vert
T_{n}y_{n}-T_{n}x_{n}\right\Vert ^{2}  \label{oo} \\
&\leq &\alpha _{n}\left\Vert T_{n}y_{n}-T_{n}x^{\ast }\right\Vert
^{2}+\left( 1-\alpha _{n}\right) \left\Vert T_{n}x_{n}-T_{n}x^{\ast
}\right\Vert ^{2}  \notag \\
&\leq &\alpha _{n}\left\Vert y_{n}-x^{\ast }\right\Vert ^{2}+\left( 1-\alpha
_{n}\right) \left\Vert x_{n}-x^{\ast }\right\Vert ^{2}  \notag \\
&\leq &\alpha _{n}\left\Vert z_{n}-x^{\ast }\right\Vert ^{2}+\left( 1-\alpha
_{n}\right) \left\Vert x_{n}-x^{\ast }\right\Vert ^{2}  \notag \\
&\leq &\alpha _{n}\left[ \left\Vert x_{n}-x^{\ast }\right\Vert ^{2}+2\theta
_{n}\left\Vert x_{n}-x^{\ast }\right\Vert \left\Vert
x_{n}-x_{n-1}\right\Vert \right.  \notag \\
&&\left. +\theta _{n}^{2}\left\Vert x_{n}-x_{n-1}\right\Vert ^{2}\right]
+\left( 1-\alpha _{n}\right) \left\Vert x_{n}-x^{\ast }\right\Vert ^{2}
\notag \\
&=&\left\Vert x_{n}-x^{\ast }\right\Vert ^{2}+2\alpha _{n}\theta
_{n}\left\Vert x_{n}-x^{\ast }\right\Vert \left\Vert x_{n}-x_{n-1}\right\Vert
\label{6} \\
&&+\alpha _{n}\theta _{n}^{2}\left\Vert x_{n}-x_{n-1}\right\Vert ^{2}.
\notag
\end{eqnarray}

Also again by using the properties $(2)$ and $\left( 3\right) $ of Lemma (%
\ref{a2}), we get
\begin{eqnarray*}
\left\Vert x_{n+1}-x^{\ast }\right\Vert ^{2} &=&\left\Vert \left( 1-\gamma
_{n}\right) T_{n}\omega _{n}+\gamma _{n}f\left( \omega _{n}\right) -x^{\ast
}\right\Vert ^{2} \\
&\leq &\left\Vert \left( 1-\gamma _{n}\right) \left( T_{n}\omega
_{n}-x^{\ast }\right) +\gamma _{n}\left( f\left( \omega _{n}\right) -f\left(
x^{\ast }\right) \right) +\gamma _{n}\left( f\left( x^{\ast }\right)
-x^{\ast }\right) \right\Vert ^{2} \\
&\leq &\left\Vert \left( 1-\gamma _{n}\right) \left( T_{n}\omega
_{n}-x^{\ast }\right) +\gamma _{n}\left( f\left( \omega _{n}\right) -f\left(
x^{\ast }\right) \right) \right\Vert ^{2} \\
&&+2\left\langle \gamma _{n}\left( f\left( x^{\ast }\right) -x^{\ast
}\right) ,x_{n+1}-x^{\ast }\right\rangle \\
&\leq &\left( 1-\gamma _{n}\right) \left\Vert \omega _{n}-x^{\ast
}\right\Vert ^{2}+\gamma _{n}\left\Vert f\left( \omega _{n}\right) -f\left(
x^{\ast }\right) \right\Vert ^{2} \\
&&-\gamma _{n}\left( 1-\gamma _{n}\right) \left\Vert f\left( \omega
_{n}\right) -f\left( x^{\ast }\right) -\omega _{n}-x^{\ast }\right\Vert \\
&&+2\left\langle \gamma _{n}\left( f\left( x^{\ast }\right) -x^{\ast
}\right) ,x_{n+1}-x^{\ast }\right\rangle \\
&\leq &\left( 1-\gamma _{n}\right) \left\Vert \omega _{n}-x^{\ast
}\right\Vert ^{2}+\gamma _{n}k^{2}\left\Vert \omega _{n}-x^{\ast
}\right\Vert ^{2} \\
&&+2\gamma _{n}\left\langle f\left( x^{\ast }\right) -x^{\ast
},x_{n+1}-x^{\ast }\right\rangle .
\end{eqnarray*}

So, It follows from (\ref{6}) that,%
\begin{eqnarray}
\left\Vert x_{n+1}-x^{\ast }\right\Vert ^{2} &\leq &\left( 1-\gamma
_{n}\left( 1-k\right) \right) \left\Vert \omega _{n}-x^{\ast }\right\Vert
^{2}+2\gamma _{n}\left\langle f\left( x^{\ast }\right) -x^{\ast
},x_{n+1}-x^{\ast }\right\rangle  \notag \\
&=&\left( 1-\gamma _{n}\left( 1-k\right) \right) \left[ \left\Vert
x_{n}-x^{\ast }\right\Vert ^{2}+2\alpha _{n}\theta _{n}\left\Vert
x_{n}-x^{\ast }\right\Vert \left\Vert x_{n}-x_{n-1}\right\Vert \right.
\notag \\
&&\left. +\alpha _{n}\theta _{n}^{2}\left\Vert x_{n}-x_{n-1}\right\Vert ^{2}
\right] +2\gamma _{n}\left\langle f\left( x^{\ast }\right) -x^{\ast
},x_{n+1}-x^{\ast }\right\rangle  \notag \\
&\leq &\left( 1-\gamma _{n}\left( 1-k\right) \right) \left\Vert
x_{n}-x^{\ast }\right\Vert ^{2}+\alpha _{n}\theta _{n}\left\Vert
x_{n}-x_{n-1}\right\Vert \left[ 2\left\Vert x_{n}-x^{\ast }\right\Vert
+\theta _{n}\left\Vert x_{n}-x_{n-1}\right\Vert \right]  \notag \\
&&+2\gamma _{n}\left\langle f\left( x^{\ast }\right) -x^{\ast
},x_{n+1}-x^{\ast }\right\rangle .  \label{7}
\end{eqnarray}

Since
\begin{equation}
\theta _{n}\left\Vert x_{n}-x_{n-1}\right\Vert =\gamma _{n}.\frac{\theta _{n}%
}{\gamma _{n}}\left\Vert x_{n}-x_{n-1}\right\Vert \rightarrow 0\text{ as }%
n\rightarrow \infty ,  \label{17}
\end{equation}

there exists a positive constant $M_{2}>0$ such that
\begin{equation*}
\theta _{n}\left\Vert x_{n}-x_{n-1}\right\Vert \leq M_{2},
\end{equation*}%
for all $n\geq 1.$From (\ref{7}), we can write%
\begin{eqnarray*}
\left\Vert x_{n+1}-x^{\ast }\right\Vert ^{2} &\leq &\left( 1-\gamma
_{n}\left( 1-k\right) \right) \left\Vert x_{n}-x^{\ast }\right\Vert
^{2}+3M_{3}a^{^{\prime }}\theta _{n}\left\Vert x_{n}-x_{n-1}\right\Vert \\
&&+2\gamma _{n}\left\langle f\left( x^{\ast }\right) -x^{\ast
},x_{n+1}-x^{\ast }\right\rangle \\
&\leq &\left( 1-\gamma _{n}\left( 1-k\right) \right) \left\Vert
x_{n}-x^{\ast }\right\Vert ^{2}+\gamma _{n}\left( 1-k\right) \left[ \frac{%
3M_{3}a^{^{\prime }}}{\left( 1-k\right) }.\frac{\theta _{n}}{\gamma _{n}}%
\left\Vert x_{n}-x_{n-1}\right\Vert \right. \\
&&\left. \frac{2}{\left( 1-k\right) }2\gamma _{n}\left\langle f\left(
x^{\ast }\right) -x^{\ast },x_{n+1}-x^{\ast }\right\rangle \right] ,
\end{eqnarray*}

where $M_{3}=\sup_{n\geq 1}\left\{ \left\Vert x_{n}-x^{\ast }\right\Vert
,M_{2}\right\} .$ In the above inequality, if we set
\begin{equation*}
s_{n}=\left\Vert x_{n}-x^{\ast }\right\Vert ^{2},\text{ }\delta _{n}=\gamma
_{n}\left( 1-k\right)
\end{equation*}

and
\begin{equation*}
r_{n}=\frac{3M_{3}a^{^{\prime }}}{\left( 1-k\right) }\frac{\theta _{n}}{%
\gamma _{n}}\left\Vert x_{n}-x_{n-1}\right\Vert +\frac{2}{\left( 1-k\right) }%
2\gamma _{n}\left\langle f\left( x^{\ast }\right) -x^{\ast },x_{n+1}-x^{\ast
}\right\rangle ,
\end{equation*}

then we obtain%
\begin{equation}
s_{n+1}\leq \left( 1-\delta _{n}\right) s_{n}+\delta _{n}r_{n},\forall n\geq
1.\text{ }  \label{18}
\end{equation}

Now we need to prove $\limsup_{n\rightarrow \infty }r_{n}\leq 0$\ in order
to complete the proof. So, we consider the following two cases.

In the first case, we assume that there exists $n_{0}\in
\mathbb{N}
$ such that the sequence $\left\{ \left\Vert x_{n}-x^{\ast }\right\Vert
\right\} _{n\geq n_{0}}$ is nonincreasing. Since the sequence $\left\{
x_{n}\right\} $ is bounded, it follows that $\left\{ \left\Vert
x_{n}-x^{\ast }\right\Vert \right\} $ is a convergent sequence. By using the
assumption $\left( 3\right) $ of theorem, we get $\sum\limits_{n=1}^{\infty
}\delta _{n}=\infty .$ Next, we claim that%
\begin{equation*}
\limsup_{n\rightarrow \infty }\left\langle f\left( x^{\ast }\right) -x^{\ast
},x_{n+1}-x^{\ast }\right\rangle \leq 0.
\end{equation*}

Since%
\begin{eqnarray}
\left\Vert y_{n}-x^{\ast }\right\Vert ^{2} &=&\left\Vert \left( 1-\beta
_{n}\right) z_{n}+\beta _{n}T_{n}z_{n}-x^{\ast }\right\Vert ^{2}  \notag \\
&=&\left\Vert \left( 1-\beta _{n}\right) \left( z_{n}-x^{\ast }\right)
+\beta _{n}\left( T_{n}z_{n}-x^{\ast }\right) \right\Vert ^{2}  \notag \\
&\leq &\beta _{n}\left\Vert T_{n}z_{n}-x^{\ast }\right\Vert ^{2}+\left(
1-\beta _{n}\right) \left\Vert z_{n}-x^{\ast }\right\Vert ^{2}-\beta
_{n}\left( 1-\beta _{n}\right) \left\Vert T_{n}z_{n}-z_{n}\right\Vert ^{2}
\notag \\
&\leq &\left\Vert z_{n}-x^{\ast }\right\Vert ^{2}-\beta _{n}\left( 1-\beta
_{n}\right) \left\Vert T_{n}z_{n}-z_{n}\right\Vert ^{2},  \label{8}
\end{eqnarray}

it follows from (\ref{5}) and the property $\left( 4\right) $ of Lemma (\ref%
{a2}) that
\begin{eqnarray}
\left\Vert \omega _{n}-x^{\ast }\right\Vert ^{2} &\leq &\alpha
_{n}\left\Vert T_{n}y_{n}-T_{n}x^{\ast }\right\Vert ^{2}+\left( 1-\alpha
_{n}\right) \left\Vert T_{n}x_{n}-T_{n}x^{\ast }\right\Vert ^{2}  \notag \\
&&-\alpha _{n}\left( 1-\alpha _{n}\right) \left\Vert
T_{n}y_{n}-T_{n}x_{n}\right\Vert ^{2}  \notag \\
&\leq &\alpha _{n}\left\Vert z_{n}-x^{\ast }\right\Vert ^{2}-\alpha
_{n}\beta _{n}\left( 1-\beta _{n}\right) \left\Vert
T_{n}z_{n}-z_{n}\right\Vert ^{2}+\left( 1-\alpha _{n}\right) \left\Vert
x_{n}-x^{\ast }\right\Vert ^{2}  \notag \\
&\leq &\alpha _{n}\left\Vert x_{n}-x^{\ast }\right\Vert ^{2}+2\alpha
_{n}\theta _{n}\left\Vert x_{n}-x^{\ast }\right\Vert \left\Vert
x_{n}-x_{n-1}\right\Vert +\alpha _{n}\theta _{n}^{2}\left\Vert
x_{n}-x_{n-1}\right\Vert ^{2}  \notag \\
&&-\alpha _{n}\beta _{n}\left( 1-\beta _{n}\right) \left\Vert
T_{n}z_{n}-z_{n}\right\Vert ^{2}+\left( 1-\alpha _{n}\right) \left\Vert
x_{n}-x^{\ast }\right\Vert ^{2}  \notag \\
&=&\left\Vert x_{n}-x^{\ast }\right\Vert ^{2}+2\alpha _{n}\theta
_{n}\left\Vert x_{n}-x^{\ast }\right\Vert \left\Vert x_{n}-x_{n-1}\right\Vert
\notag \\
&&+\alpha _{n}\theta _{n}^{2}\left\Vert x_{n}-x_{n-1}\right\Vert ^{2}-\alpha
_{n}\beta _{n}\left( 1-\beta _{n}\right) \left\Vert
T_{n}z_{n}-z_{n}\right\Vert ^{2}.  \label{9'}
\end{eqnarray}

Also, from (\ref{9'}) and the property $(3)$ of Lemma (\ref{a2}), we get%
\begin{eqnarray*}
\left\Vert x_{n+1}-x^{\ast }\right\Vert ^{2} &=&\left\Vert \left( 1-\gamma
_{n}\right) \left( T_{n}\omega _{n}-x^{\ast }\right) +\gamma _{n}\left(
f\left( \omega _{n}\right) -x^{\ast }\right) \right\Vert ^{2} \\
&\leq &\gamma _{n}\left\Vert f\left( \omega _{n}\right) -x^{\ast
}\right\Vert ^{2}+\left( 1-\gamma _{n}\right) \left\Vert T_{n}\omega
_{n}-T_{n}x^{\ast }\right\Vert ^{2} \\
&&-\gamma _{n}\left( 1-\gamma _{n}\right) \left\Vert f\left( \omega
_{n}\right) -T_{n}\omega _{n}\right\Vert ^{2} \\
&\leq &\gamma _{n}\left\Vert f\left( \omega _{n}\right) -x^{\ast
}\right\Vert ^{2}+\left( 1-\gamma _{n}\right) \left\Vert \omega _{n}-x^{\ast
}\right\Vert ^{2} \\
&\leq &\gamma _{n}\left\Vert f\left( \omega _{n}\right) -x^{\ast
}\right\Vert ^{2}+\left( 1-\gamma _{n}\right) \left[ \left\Vert
x_{n}-x^{\ast }\right\Vert ^{2}+2\alpha _{n}\theta _{n}\left\Vert
x_{n}-x^{\ast }\right\Vert \left\Vert x_{n}-x_{n-1}\right\Vert \right. \\
&&\left. +\alpha _{n}\theta _{n}^{2}\left\Vert x_{n}-x_{n-1}\right\Vert
^{2}-\alpha _{n}\beta _{n}\left( 1-\beta _{n}\right) \left\Vert
T_{n}z_{n}-z_{n}\right\Vert ^{2}\right] \\
&=&\gamma _{n}\left\Vert f\left( \omega _{n}\right) -x^{\ast }\right\Vert
^{2}+\left( 1-\gamma _{n}\right) \left\Vert x_{n}-x^{\ast }\right\Vert
^{2}+\left( 1-\gamma _{n}\right) \alpha _{n}\theta _{n}\left\Vert
x_{n}-x_{n-1}\right\Vert \left[ 2\left\Vert x_{n}-x^{\ast }\right\Vert
\right. \\
&&\left. \theta _{n}\left\Vert x_{n}-x_{n-1}\right\Vert \right] -\left(
1-\gamma _{n}\right) \alpha _{n}\beta _{n}\left( 1-\beta _{n}\right)
\left\Vert T_{n}z_{n}-z_{n}\right\Vert ^{2}.
\end{eqnarray*}

So the following is true for all $n>1:$
\begin{eqnarray}
\left( 1-\gamma _{n}\right) \alpha _{n}\beta _{n}\left( 1-\beta _{n}\right)
\left\Vert T_{n}z_{n}-z_{n}\right\Vert ^{2} &\leq &\gamma _{n}\left(
\left\Vert f\left( \omega _{n}\right) -x^{\ast }\right\Vert ^{2}-\left\Vert
x_{n}-x^{\ast }\right\Vert ^{2}\right) +\left\Vert x_{n}-x^{\ast
}\right\Vert ^{2}-\left\Vert x_{n+1}-x^{\ast }\right\Vert ^{2}  \notag \\
&&\left( 1-\gamma _{n}\right) \alpha _{n}\theta _{n}\left\Vert
x_{n}-x_{n-1}\right\Vert \left[ 2\theta _{n}\left\Vert x_{n}-x^{\ast
}\right\Vert +\theta _{n}\left\Vert x_{n}-x_{n-1}\right\Vert \right] .
\notag
\end{eqnarray}

Hence, we conclude from the assumptions $\left( 1\right) $ and $\left(
2\right) $ of theorem and the convergence of the sequences $\left\{
\left\Vert x_{n}-x^{\ast }\right\Vert \right\} $ and of $\left\{ \theta
_{n}\left\Vert x_{n}-x_{n-1}\right\Vert \right\} $ that%
\begin{equation}
\left\Vert T_{n}z_{n}-z_{n}\right\Vert \rightarrow 0\text{ as }n\rightarrow
\infty .  \label{11}
\end{equation}

Since $\left\{ T_{n}\right\} $ satisfies the NST-condition (I) with $T,$ we
obtain that%
\begin{equation}
\left\Vert Tz_{n}-z_{n}\right\Vert \rightarrow 0\text{ as }n\rightarrow
\infty \text{ }.  \label{11a}
\end{equation}

On the other hand, since%
\begin{equation}
\left\Vert z_{n}-x_{n}\right\Vert =\left\Vert x_{n}+\theta _{n}\left(
x_{n}-x_{n-1}\right) -x_{n}\right\Vert =\theta _{n}\left\Vert
x_{n}-x_{n-1}\right\Vert \rightarrow 0\text{ as }n\rightarrow \infty ,
\label{12}
\end{equation}

by using (\ref{11}) and (\ref{12}), we have%
\begin{eqnarray*}
\left\Vert y_{n}-x_{n}\right\Vert &=&\left\Vert \left( 1-\beta _{n}\right)
z_{n}+\beta _{n}T_{n}z_{n}-x_{n}\right\Vert \\
&\leq &\left\Vert z_{n}-x_{n}\right\Vert +\beta _{n}\left\Vert
Tz_{n}-z_{n}\right\Vert ,
\end{eqnarray*}%
and so it is provided that
\begin{equation}
\left\Vert y_{n}-x_{n}\right\Vert \rightarrow 0\text{ as }n\rightarrow
\infty .  \label{14}
\end{equation}

Also, from (\ref{11}), (\ref{12}) and (\ref{14}), we get
\begin{eqnarray*}
\left\Vert \omega _{n}-z_{n}\right\Vert &=&\left\Vert \omega
_{n}-z_{n}+T_{n}z_{n}-T_{n}z_{n}\right\Vert \\
&\leq &\left\Vert \omega _{n}-T_{n}z_{n}\right\Vert +\left\Vert
T_{n}z_{n}-z_{n}\right\Vert \\
&=&\left\Vert T_{n}\left( \left( 1-\alpha _{n}\right) T_{n}x_{n}+\alpha
_{n}T_{n}y_{n}\right) -T_{n}z_{n}\right\Vert +\left\Vert
T_{n}z_{n}-z_{n}\right\Vert \\
&\leq &\left\Vert \left( 1-\alpha _{n}\right) T_{n}x_{n}+\alpha
_{n}T_{n}y_{n}-z_{n}\right\Vert +\left\Vert T_{n}z_{n}-z_{n}\right\Vert \\
&\leq &\alpha _{n}\left\Vert T_{n}y_{n}-T_{n}x_{n}\right\Vert +\left\Vert
T_{n}x_{n}-T_{n}z_{n}\right\Vert +\left\Vert T_{n}z_{n}-z_{n}\right\Vert
+\left\Vert T_{n}z_{n}-z_{n}\right\Vert \\
&\leq &\alpha _{n}\left\Vert y_{n}-x_{n}\right\Vert +\left\Vert
x_{n}-z_{n}\right\Vert +2\left\Vert T_{n}z_{n}-z_{n}\right\Vert ,
\end{eqnarray*}%
and so%
\begin{equation}
\left\Vert \omega _{n}-z_{n}\right\Vert \rightarrow 0\text{ as }n\rightarrow
\infty .  \label{14'}
\end{equation}%
By using (\ref{11}) and (\ref{14'}), we obtain%
\begin{eqnarray*}
\left\Vert x_{n+1}-z_{n}\right\Vert &\leq &\left\Vert
x_{n+1}-T_{n}z_{n}\right\Vert +\left\Vert z_{n}-T_{n}z_{n}\right\Vert \\
&=&\left\Vert \left( 1-\gamma _{n}\right) T_{n}\omega _{n}+\gamma
_{n}f\left( \omega _{n}\right) -T_{n}z_{n}\right\Vert +\left\Vert
T_{n}z_{n}-z_{n}\right\Vert \\
&\leq &\gamma _{n}\left\Vert f\left( \omega _{n}\right) -T_{n}\omega
_{n}\right\Vert +\left\Vert T_{n}\omega _{n}-T_{n}z_{n}\right\Vert
+\left\Vert T_{n}z_{n}-z_{n}\right\Vert \\
&\leq &\gamma _{n}\left\Vert f\left( \omega _{n}\right) -T_{n}\omega
_{n}\right\Vert +\left\Vert \omega _{n}-z_{n}\right\Vert +\left\Vert
T_{n}z_{n}-z_{n}\right\Vert
\end{eqnarray*}

which implies%
\begin{equation}
\left\Vert x_{n+1}-z_{n}\right\Vert \rightarrow 0\text{ as }n\rightarrow
\infty .  \label{13}
\end{equation}

Hence, we get%
\begin{equation*}
\left\Vert x_{n+1}-x_{n}\right\Vert \leq \left\Vert x_{n+1}-z_{n}\right\Vert
+\left\Vert z_{n}-x_{n}\right\Vert \rightarrow 0\text{ as }n\rightarrow
\infty .
\end{equation*}

Now, we are in a position to prove that $\left\{ x_{n}\right\} \rightarrow
x^{\ast }.$ Let%
\begin{equation*}
\upsilon =\limsup_{n\rightarrow \infty }\left\langle f\left( x^{\ast
}\right) -x^{\ast },x_{n+1}-x^{\ast }\right\rangle .
\end{equation*}%
Since $\left\{ x_{n}\right\} $ is bounded, there exists a subsequence $%
\left\{ x_{n_{i}}\right\} $\ of $\left\{ x_{n}\right\} $\ such that $%
x_{n_{i}}\rightharpoonup t$ and \textbf{\ }%
\begin{equation*}
\upsilon =\lim_{i\rightarrow \infty }\left\langle f\left( x^{\ast }\right)
-x^{\ast },x_{n_{i}+1}-x^{\ast }\right\rangle .
\end{equation*}

By using (\ref{11a}) and (\ref{12}), we obtain that%
\begin{eqnarray*}
\left\Vert x_{n}-Tx_{n}\right\Vert &=&\left\Vert
x_{n}-z_{n}+z_{n}+Tz_{n}-Tz_{n}-Tx_{n}\right\Vert \\
&\leq &\left\Vert z_{n}-x_{n}\right\Vert +\left\Vert Tz_{n}-z_{n}\right\Vert
+\left\Vert Tz_{n}-Tx_{n}\right\Vert \\
&\leq &2\left\Vert z_{n}-x_{n}\right\Vert +\left\Vert Tz_{n}-z_{n}\right\Vert
\end{eqnarray*}

and so, we conclude that%
\begin{equation*}
\left\Vert x_{n}-Tx_{n}\right\Vert \rightarrow 0\text{ as }n\rightarrow
\infty .
\end{equation*}

Hence, it is obvious from Lemma \ref{lemmaaa} that $t\in F\left( T\right) .$
On the other hand, since $\left\Vert x_{n+1}-x_{n}\right\Vert \rightarrow 0$%
\ as $n\rightarrow \infty $\ and $x_{n_{i}}\rightharpoonup t$\ this implies
that $x_{n_{i}+1}\rightarrow t.$\textbf{\ }Furthermore, by using $x^{\ast
}=P_{F\left( T\right) }f\left( x^{\ast }\right) $ and the property of the
metric projection operators (\ref{30}), we can \textbf{w}rite%
\begin{equation}
\upsilon =\lim_{i\rightarrow \infty }\left\langle f\left( x^{\ast }\right)
-x^{\ast },x_{n_{i}+1}-x^{\ast }\right\rangle =\left\langle f\left( x^{\ast
}\right) -x^{\ast },t-x^{\ast }\right\rangle \leq 0.  \label{15}
\end{equation}

Then, we have%
\begin{equation}
\limsup_{n\rightarrow \infty }\left\langle f\left( x^{\ast }\right) -x^{\ast
},x_{n+1}-x^{\ast }\right\rangle \leq 0.  \label{16}
\end{equation}

It follows from (\ref{17}) and (\ref{16}) that $\limsup_{n\rightarrow \infty
}r_{n}\leq 0.$ Finally, we deduce that $x_{n}\rightarrow x^{\ast }.$

In the second case, we assume that there exists a $n_{0}\in
\mathbb{N}
$ such that the sequence $\left\{ \left\Vert x_{n}-x^{\ast }\right\Vert
\right\} _{n\geq n_{0}}$ is not monotonically decreasing. Let $\Phi
_{n}=\left\Vert x_{n}-x^{\ast }\right\Vert ^{2}$ for all $n\in
\mathbb{N}
$. So, there exists a subsequence $\left\{ \Phi _{n_{j}}\right\} $ of $%
\left\{ \Phi _{n}\right\} $ such that $\Phi _{n_{j}}<\Phi _{n_{j+1}}$ for
all $j\in
\mathbb{N}
.$ Then, we define $\tau :\left\{ n:n\geq n_{0}\right\} \rightarrow
\mathbb{N}
$ as follows:
\begin{equation*}
\tau \left( n\right) :=\max \left\{ l\in
\mathbb{N}
:l\leq n,\Phi _{l}<\Phi _{l+1}\right\} .
\end{equation*}

Obviously, $\tau $ is a nondecreasing sequence. Then, by using Lemma \ref%
{main}, we have $\Phi _{\tau \left( n\right) }\leq \Phi _{\tau \left(
n\right) +1}$, i.e., $\left\Vert x_{\tau \left( n\right) }-x^{\ast
}\right\Vert \leq \left\Vert x_{\tau \left( n\right) +1}-x^{\ast
}\right\Vert $ for all $n\geq n_{0}.$ Similarly to the first case, we can
obtain the everything proved in the first case by taking $\tau \left(
n\right) $ instead of $n$. So we have
\begin{equation*}
\limsup_{n\rightarrow \infty }\left\Vert x_{\tau \left( n\right) }-x^{\ast
}\right\Vert ^{2}\leq 0.
\end{equation*}

Therefore, we get%
\begin{equation}
\left\Vert x_{\tau \left( n\right) }-x^{\ast }\right\Vert ^{2}\rightarrow 0\
\text{\ and }\left\Vert x_{\tau \left( n\right) +1}-x^{\ast }\right\Vert
\rightarrow 0\text{ as }n\rightarrow \infty .  \label{22}
\end{equation}

So, it follows from (\ref{22}) and Lemma \ref{main} that%
\begin{equation*}
\left\Vert x_{n}-x^{\ast }\right\Vert \leq \left\Vert x_{\tau \left(
n\right) +1}-x^{\ast }\right\Vert \rightarrow 0\text{ as }n\rightarrow
\infty .
\end{equation*}

Therefore, $\left\{ x_{n}\right\} $ converges strongly to $x^{\ast }$ and
this completes the proof.
\end{proof}

\section{Applications}

In this section,\ we will give some applications of Algorithm \ref{1''} to
the convex minimization, variational inequality, generalized equilibrium,
monotone inclusion and image restoration problems.

\subsection{Applicaton to monotone inclusion problems}

Finding the zero sum of two monotone operators is one of the most important
problems in monotone operator theory. We study the following inclusion
problem: finding $x\in H$ such that%
\begin{equation}
0\in \left( A+B\right) x  \label{1'}
\end{equation}%
where $A:H\rightarrow H$ is an operator and $B:H\rightarrow 2^{H}$ is a
set-valued operator. This problem includes, as special cases, convex
minimization problems, variational inequalities, and equilibrium problems.
Also, some tangible problems in statistical regression, machine learning,
image processing, signal processing and the linear inverse problem can be
formulated mathematically in the form (\ref{1'}). It is well-known that the
problem (\ref{1'}) is equivalent to the problem of finding $x$ which
satisfies the following equation:%
\begin{equation*}
J_{\lambda }^{B}\left( I-\lambda A\right) x=x,
\end{equation*}%
where $A:H\rightarrow H$ a $\eta $-inverse strongly monotone operator, $%
B:H\rightarrow 2^{H}$ a maximal monotone operator and $\lambda \in \left[
0,2\eta \right] .$ Also, it can be seen that $J_{\lambda }^{B}\left(
I-\lambda A\right) $ is a nonexpansive mapping, see for details \cite{bau}.

Now, as a corollary of Theorem \ref{theorem}, we give the following to
approximate a solution of the inclusion problem (\ref{1'}) by swaping $T_{n}$
and $T$ with $J_{\lambda _{n}}^{B}\left( I-\lambda _{n}A\right) $ and $%
J_{\lambda }^{B}\left( I-\lambda A\right) ,$ respectively.

\begin{theorem}
\label{coral} Let $f$ be a $k$-contraction mapping on $H$, $A:H\rightarrow H$
a $\eta $-inverse strongly monotone operator and $B:H\rightarrow 2^{H}$ a
maximal monotone operator such that $\Omega =\left( A+B\right) ^{-1}\left(
0\right) \neq \emptyset .$ Let $\left\{ \lambda _{n}\right\} \in \left(
0,2\eta \right) $ be a sequence such that $\lambda _{n}\rightarrow \lambda $
where $\lambda $ is a constant belongs to $\left( 0,2\eta \right) .$ Let $%
x_{0},$ $x_{1}\in H,$ $\theta \geq 0$ and $\left\{ x_{n}\right\} $ be a
sequence generated by%
\begin{equation}
\left\{
\begin{array}{l}
z_{n}=x_{n}+\theta _{n}\left( x_{n}-x_{n-1}\right) \\
y_{n}=\left( 1-\beta _{n}\right) z_{n}+\beta _{n}J_{\lambda _{n}}^{B}\left(
I-\lambda _{n}A\right) z_{n} \\
\omega _{n}=J_{\lambda _{n}}^{B}\left( I-\lambda _{n}A\right) \left( \left(
1-\alpha _{n}\right) J_{\lambda _{n}}^{B}\left( I-\lambda _{n}A\right)
x_{n}+\alpha _{n}J_{\lambda _{n}}^{B}\left( I-\lambda _{n}A\right)
y_{n}\right) \\
x_{n+1}=\left( 1-\gamma _{n}\right) J_{\lambda _{n}}^{B}\left( I-\lambda
_{n}A\right) \omega _{n}+\gamma _{n}f\left( \omega _{n}\right) .%
\end{array}%
\right.  \label{200}
\end{equation}%
where all the parameters satisfy the same conditions as in Theorem \ref%
{theorem} . Then $\left\{ x_{n}\right\} $ converges strongly to a point $%
x^{\ast }$ of $\ \Omega $, where $x^{\ast }=P_{\Omega }f\left( x^{\ast
}\right) .$
\end{theorem}

\begin{proof}
The proof is clear from Lemma \ref{FORWARD}.
\end{proof}

We next give a numerical example to show the convergence of the sequence
generated by Algorithm \ref{200} to the solution of variational inclusion
problem.

\begin{figure}[tbph]
\begin{center}
\includegraphics[width=9.0cm]{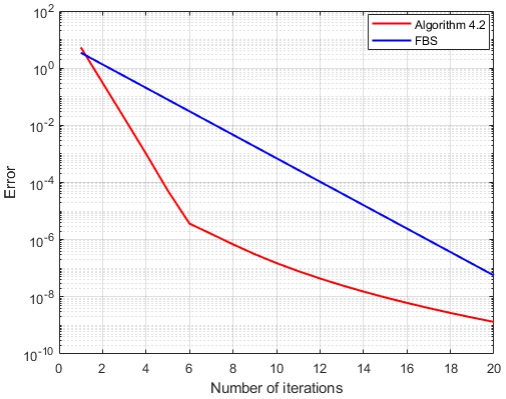}
\end{center}
\caption{ Comparison of Algorithm \protect\ref{200} and FBS algorithm}
\label{fig2b}
\end{figure}

\begin{example}
\label{exam} Let $A:l_{2}\rightarrow l_{2}$ and $B:l_{2}\rightarrow l_{2}$
be two operators defined by $Ax=3x+\left( 1,2,3,0,...\right) $ and $Bx=8x,$
where $x=\left( x_{1},x_{2},x_{3,}...\right) \in l_{2}.$ We can easily see
that $A$ is a $1/3$-inverse strongly monotone and $B$ is a maximal monotone
operator. Indeed, for $x,y\in l_{2},$ we have%
\begin{equation*}
\left\langle x-y,Ax-Ay\right\rangle =\left\langle x-y,3x-3y\right\rangle
=3\left\Vert x-y\right\Vert _{l_{2}}^{2}\geq \frac{1}{3}\left\Vert
Ax-Ay\right\Vert _{l_{2}}^{2}
\end{equation*}

and
\begin{equation*}
\left\langle x-y,Bx-By\right\rangle =\left\langle x-y,8x-8y\right\rangle
=8\left\Vert x-y\right\Vert _{l_{2}}^{2}\geq 0.
\end{equation*}

So $B$ is a monotone operator. On the other hand, since $R\left( I+\lambda
_{n}B\right) =l_{2},$ we obtain from Lemma \ref{minty} that $B$ is a maximal
monotone operator. By a direct calculation, we obtain%
\begin{eqnarray*}
J_{\lambda _{n}}^{B}\left( x-\lambda _{n}Ax\right) &=&\left( I+\lambda
A\right) ^{-1}(x-\lambda _{n}Ax) \\
&=&\frac{1-3\lambda _{n}}{1+8\lambda _{n}}x-\frac{\lambda _{n}}{1+8\lambda
_{n}}\left( 1,2,3,0,...\right) ,
\end{eqnarray*}

and $\left( A+B\right) ^{-1}\left( 0\right) =(-0.0909,-0.1818,-0.2727,0,...)$%
. In Algorithm \ref{200}, we choose $x_{1}=\left( -3,-5,-1,0,...\right)
,\alpha _{n}=\frac{1}{10^{8}},$ $\beta _{n}=\frac{1}{n+1},$ $\gamma _{n}=%
\frac{1}{\left( n+1\right) ^{6}},$ $\eta _{n}=\frac{10}{n}$ and $f\left(
x\right) =0.1x.$ If $\theta =0.99,$ then we get
\begin{equation*}
\theta _{n}:=\left\{
\begin{array}{ll}
\min \left\{ 0.99,\frac{10}{n\left( n+1\right) \left\Vert
x_{n}-x_{n-1}\right\Vert }\right\} & if\text{ \ }x_{n}\neq x_{n-1}, \\
0.99, & otherwise.%
\end{array}%
\right. .
\end{equation*}%
In Table \ref{fig11b}, we give the iteration steps of Algorithm \ref{200}.

\begin{table}[tbph]
\begin{center}
\begin{tabular}{ccc}
\hline
No.Iteration & $x_{n}$ & $\left\Vert x_{n+1}-x_{n}\right\Vert $ \\ \hline
$1$ & $(-3.0000,-5.0000,-1.0000,0,...)$ & $.$ \\
$2$ & $(-0.2587,-0.4593,-0.3112,0,...)$ & $5.348508$ \\
$3$ & $(-0.1007,-0.1979,-0.2746,0,...)$ & $0.307697$ \\
$4$ & $(-0.0915,-0.1827,-0.2728,0,...)$ & $0.017842$ \\
$5$ & $(-0.0909,-0.1819,-0.2727,0,...)$ & $0.001013$ \\
$6$ & $(-0.0909,-0.1819,-0.2727,0,...)$ & $0.000052$ \\
$7$ & $(-0.0909,-0.1819,-0.2727,0,...)$ & $0.000004$ \\
$8$ & $(-0.0909,-0.1819,-0.2727,0,...)$ & $0.000002$ \\
$9$ & $(-0.0909,-0.1819,-0.2727,0,...)$ & $0.000001$ \\
$10$ & $(-0.0909,-0.1819,-0.2727,0,...)$ & $0.000000$ \\
$\vdots $ & $\vdots $ & $\vdots $ \\
$20$ & $(-0.0909,-0.1818,-0.2727,0,...)$ & $0.000000$ \\
$21$ & $(-0.0909,-0.1818,-0.2727,0,...)$ & $0.000000$ \\ \hline
\end{tabular}%
\end{center}
\caption{Iteration steps and error terms of Algorithm \protect\ref{200}}
\label{fig11b}
\end{table}
In Figure \ref{fig2b}, we compare the performance of Algorithm \ref{200} and
the FBS algorithm.
\end{example}

\subsection{Applicaton to convex minimization problems}

Let $h:H\rightarrow
\mathbb{R}
$ be a differentiable and convex function such that $\nabla h$ is a $\eta $%
-inverse strongly monotone and $g:H\rightarrow
\mathbb{R}
$ be a proper convex and lower semi-continuous function. We consider the
following convex minimization problem: finding $x^{\ast }\in H$ such that%
\begin{equation}
h\left( x^{\ast }\right) +g\left( x^{\ast }\right) =\min_{x\in H}\left\{
h\left( x\right) +g\left( x\right) \right\} .  \label{23}
\end{equation}%
The solution set of convex minimization problem is denoted by $\Pi .$ As
mentioned in the introduction part, it is known that the convex minimization
problem (\ref{23}) is equivalent to the fixed point problem (\ref{105}).
Also, since $\nabla h$ is an inverse strongly monotone and $\partial g$ is a
maximal monotone operators, we know that $J_{\lambda _{n}}^{\partial
g}\left( I-\lambda _{n}\nabla h\right) $ is a nonexpansive mapping for $%
\lambda \in \left[ 0,2\eta \right] $. So, we obtain the following theorem.

\begin{theorem}
\label{the} Let $f$ be a $k$-contraction mapping on $H,$ $h:H\rightarrow
\mathbb{R}
$ be a differentiable and convex function such that $\nabla h$ is a $\eta $%
-inverse strongly monotone and $g:H\rightarrow
\mathbb{R}
$ be a proper convex and lower semi-continuous function such that $\Pi \neq
\emptyset .$ Let $x_{0},x_{1}\in H$ and $\left\{ x_{n}\right\} $ be a
sequence generated by%
\begin{equation}
\left\{
\begin{array}{l}
z_{n}=x_{n}+\theta _{n}\left( x_{n}-x_{n-1}\right) \\
y_{n}=\left( 1-\beta _{n}\right) z_{n}+\beta _{n}J_{\lambda _{n}}^{\partial
g}\left( I-\lambda _{n}\nabla h\right) z_{n} \\
\omega _{n}=J_{\lambda _{n}}^{\partial g}\left( I-\lambda _{n}\nabla
h\right) \left( \left( 1-\alpha _{n}\right) J_{\lambda _{n}}^{\partial
g}\left( I-\lambda _{n}\nabla h\right) x_{n}+\alpha _{n}J_{\lambda
_{n}}^{\partial g}\left( I-\lambda _{n}\nabla h\right) y_{n}\right) \\
x_{n+1}=\left( 1-\gamma _{n}\right) J_{\lambda _{n}}^{\partial g}\left(
I-\lambda _{n}\nabla h\right) \omega _{n}+\gamma _{n}f\left( \omega
_{n}\right)%
\end{array}%
\right.  \label{a}
\end{equation}%
where,%
\begin{equation*}
\theta _{n}:=\left\{
\begin{array}{ll}
\min \left\{ \theta ,\frac{\eta _{n}\gamma _{n}}{\left\Vert
x_{n}-x_{n-1}\right\Vert }\right\} & if\text{ \ }x_{n}\neq x_{n-1} \\
\theta , & otherwise%
\end{array}%
\right. ,
\end{equation*}%
for $\theta \geq 0.$ Let $\left\{ \alpha _{n}\right\} ,$ $\left\{ \beta
_{n}\right\} ,$ $\left\{ \gamma _{n}\right\} ,$ $\left\{ \theta _{n}\right\}
$ and $\left\{ \eta _{n}\right\} $ be sequences which satisfy the same
conditions as in Theorem \ref{theorem} such that $\lambda _{n}\rightarrow
\lambda ,$ for $\lambda _{n},\lambda \in \left( 0,2\eta \right) .$ Then $%
\left\{ x_{n}\right\} $ converges strongly to a $x^{\ast }$ solution of
convex minimization problem, where $x^{\ast }=P_{\Pi }f\left( x^{\ast
}\right) $ $.$
\end{theorem}

\subsection{Application to the variational inequality problems (VIP)}

The variational inequality problem is defined as the problem of finding a
point $x^{\ast }\in C$ such that
\begin{equation}
\left\langle Ax^{\ast },y-x^{\ast }\right\rangle ,\text{ }\forall y\in C
\label{24}
\end{equation}%
where $A:C\rightarrow H$ is a nonlinear monotone operator. We denote the
solution set of (\ref{24}) by $VI(C,A)$. It is known that the variational
inequality problem (\ref{24}) is equivalent to finding a point $x^{\ast }$
such that $x^{\ast }=J_{\lambda _{n}}^{\partial i_{C}}\left( I-\lambda
_{n}A\right) x^{\ast }$ where $\partial i_{C}$ is the subdifferential of the
indicator function $i_{C}:H\rightarrow \left( -\infty ,\infty \right] $ of $%
C $ defined by
\begin{equation*}
i_{C}\left( x\right) =\left\{
\begin{array}{ll}
0, & \text{if }x\in C \\
\infty , & \text{if }x\notin C%
\end{array}%
\right. .
\end{equation*}%
It is well-known that the indicator function $i_{C}$ is a proper, lower
semi-continuous and convex function on $H.$ So, the subdifferential $%
\partial i_{C}$ is a maximal monotone operator. Based on these facts, we can
easily see that,%
\begin{equation*}
\begin{array}{lll}
y=J_{\lambda }^{\partial i_{C}}x & \Leftrightarrow & x\in \left( y+\lambda
\partial i_{C}y\right) \\
& \Leftrightarrow & x-y\in \lambda \partial i_{C}y \\
& \Leftrightarrow & y=P_{C}x.%
\end{array}%
\end{equation*}%
So, the variational inequality problem (\ref{24}) is equivalent to the fixed
point problem $x^{\ast }=$ $P_{C}\left( I-\lambda _{n}A\right) x^{\ast }$.
Since $P_{C}\left( I-\lambda _{n}A\right) $ is a nonexpansive mapping when $%
A $ is an inverse strongly monotone operator, the following theorem can be
obtained from Theorem \ref{theorem}.

\begin{theorem}
Let $A:H\rightarrow H$ be an $\eta $-inverse strongly monotone operator such
that $VI(C,A)\neq \emptyset .$ Let $x_{0},x_{1}\in H$ and let $\left\{
x_{n}\right\} $ be a sequence generated by%
\begin{equation*}
\left\{
\begin{array}{l}
z_{n}=x_{n}+\theta _{n}\left( x_{n}-x_{n-1}\right) \\
y_{n}=\left( 1-\beta _{n}\right) z_{n}+\beta _{n}P_{C}\left( I-\lambda
_{n}A\right) z_{n} \\
\omega _{n}=P_{C}\left( I-\lambda _{n}A\right) \left( \left( 1-\alpha
_{n}\right) P_{C}\left( I-\lambda _{n}A\right) x_{n}+\alpha _{n}P_{C}\left(
I-\lambda _{n}A\right) y_{n}\right) \\
x_{n+1}=\left( 1-\gamma _{n}\right) P_{C}\left( I-\lambda _{n}A\right)
\omega _{n}+\gamma _{n}f\left( \omega _{n}\right)%
\end{array}%
\right.
\end{equation*}%
where,
\begin{equation*}
\theta _{n}:=\left\{
\begin{array}{ll}
\min \left\{ \theta ,\frac{\eta _{n}\gamma _{n}}{\left\Vert
x_{n}-x_{n-1}\right\Vert }\right\} & if\text{ \ }x_{n}\neq x_{n-1} \\
\theta , & otherwise%
\end{array}%
\right. \text{ },
\end{equation*}%
and for $\theta \geq 0.$ Then under the same conditions as in Theorem \ref%
{theorem}, the sequence $\left\{ x_{n}\right\} $ converges strongly to a
solution $x^{\ast }$ of variational inequality problem (\ref{24}), where $%
x^{\ast }=P_{VI(C,A)}f\left( x^{\ast }\right) .$
\end{theorem}

\subsection{Application to generalized equilibrium problems (GEPs)}

Let $F:C\times C\rightarrow
\mathbb{R}
$ be a bifunction and $A:H\rightarrow H$ \ be a monotone operator. The
generalized equilibrium problem is formulated by finding $x\in C$ such that
\begin{equation}
F\left( x,y\right) +\left\langle Ax,y-x\right\rangle \geq 0  \label{25}
\end{equation}%
for all $y\in C$. In order to solve the generalized equilibrium problem (\ref%
{25}), we assume that $F$ satisfies the following:

\begin{enumerate}
\item[(1)] $F\left( x,x\right) =0,$ for all $x,y\in C,$

\item[(2)] $F$ is monotone, that is, $F\left( x,y\right) +F\left( y,x\right)
\leq 0$, for all $x,y\in C,$

\item[(3)] for each $x,y,z\in C$ $\lim_{t\rightarrow 0}F\left( tz+\left(
1-t\right) x,y\right) \leq F\left( x,y\right) ,$

\item[(4)] for each $x,y,z\in C,$ $y\mapsto F\left( x,y\right) $ is convex
and lower semicontinuous.
\end{enumerate}

We denote the set of solutions of (\ref{25}) by $GEP\left( F,A\right) $. We
need the following lemmas to give an application of our main theorem to the
generalized equilibrium problem.

\begin{lemma}
\label{uc}\cite{taka} Let $F$ be a bifunction from $C\times C$ to $%
\mathbb{R}
$ which satisfies $\left( 1\right) $-$\left( 4\right) $. Let $A_{F}$ be a
set-valued mapping from $H$ into itself defined by
\begin{equation}
A_{F}x=\left\{
\begin{array}{ll}
\left\{ z\in H:F\left( x,y\right) \geq \left\langle y-x,z\right\rangle
\right\} , & x\in C \\
\emptyset , & x\notin C.%
\end{array}%
\right.  \label{A}
\end{equation}%
Then $A_{F}$ is a maximal monotone operator with the domain $D\left(
A_{F}\right) \subset C$ and $EP\left( F\right) =A_{F}^{-1}0.$
\end{lemma}

\begin{lemma}
\cite{kit1} Let $F:C\times C\rightarrow
\mathbb{R}
$ be a bifunction satisfying $\left( 1\right) $-$\left( 2\right) $ and $%
A:H\rightarrow H$ is continuous and monotone on $H$, hence maximal monotone.
For $r>0$ and $x\in H,$ let $W_{r}:H\rightarrow C$ be a mapping defined as
follows:%
\begin{equation*}
W_{r}x=\left\{ z\in C:F\left( z,y\right) +\frac{1}{r}\left\langle
z-x+rAx,y-z\right\rangle \geq 0,\text{ }\forall y\in C\right\} .
\end{equation*}%
Then the following hold:

\begin{enumerate}
\item[(1)] $W_{r}$ is single-valued,

\item[(2)] $F\left( W_{r}\right) =GEP\left( F,A\right) ,$

\item[(3)] $W_{r}$ is a nonexpansive mapping$,$

\item[(4)] $GEP\left( F,A\right) $ is closed and convex.
\end{enumerate}
\end{lemma}

Kitkuan et al. \cite{kit1} showed in the proof of the lemma that $GEP\left(
F,A\right) =\left( A+A_{F}\right) ^{-1}\left( 0\right) $. On the other hand,
it is know that $J_{\lambda _{n}}^{A_{F}}\left( I-\lambda _{n}A\right)
x=x\Leftrightarrow x\in \left( A+A_{F}\right) ^{-1}\left( 0\right) $. So, by
taking $B=A_{F}$ in Theorem \ref{coral} and by using Lemma \ref{uc}, we
obtain the following theorem.

\begin{theorem}
Let $A:H\rightarrow H$ be an $\eta $-inverse strongly monotone operator for $%
\eta >0$, $F$ a bifunction from $C\times C$ to $%
\mathbb{R}
$ which satisfies $\left( 1\right) $-$\left( 4\right) ,$ and $%
A_{F}:H\rightarrow 2^{H}$ a maximal monotone operator defined by (\ref{A})
such that $\Gamma =\left( A+A_{F}\right) ^{-1}\left( 0\right) \neq \emptyset
.$ Let $x_{0},x_{1}\in H$ and let $\left\{ x_{n}\right\} $ be a sequence
generated by%
\begin{equation*}
\left\{
\begin{array}{l}
z_{n}=x_{n}+\theta _{n}\left( x_{n}-x_{n-1}\right) \\
y_{n}=\left( 1-\beta _{n}\right) z_{n}+\beta _{n}J_{\lambda
_{n}}^{A_{F}}\left( I-\lambda _{n}A\right) z_{n} \\
\omega _{n}=J_{\lambda _{n}}^{A_{F}}\left( I-\lambda _{n}A\right) \left(
\left( 1-\alpha _{n}\right) J_{\lambda _{n}}^{A_{F}}\left( I-\lambda
_{n}A\right) x_{n}+\alpha _{n}J_{\lambda _{n}}^{A_{F}}\left( I-\lambda
_{n}A\right) y_{n}\right) \\
x_{n+1}=\left( 1-\gamma _{n}\right) J_{\lambda _{n}}^{A_{F}}\left( I-\lambda
_{n}A\right) \omega _{n}+\gamma _{n}f\left( \omega _{n}\right) .%
\end{array}%
\right.
\end{equation*}%
where,
\begin{equation*}
\theta _{n}:=\left\{
\begin{array}{ll}
\min \left\{ \theta ,\frac{\eta _{n}\gamma _{n}}{\left\Vert
x_{n}-x_{n-1}\right\Vert }\right\} & if\text{ \ }x_{n}\neq x_{n-1}, \\
\theta , & otherwise.%
\end{array}%
\right. \text{ }
\end{equation*}%
and for $\theta \geq 0.$Then under the same conditions as in Theorem \ref%
{theorem}, the sequence $\left\{ x_{n}\right\} $ converges strongly to a
point $x^{\ast }$of $\ $generalized equilibrium problem, where $x^{\ast
}=P_{\Gamma }f\left( x^{\ast }\right) .$
\end{theorem}

\subsection{Application to image restoration problems}

In this section, we apply Algorithm \ref{a} to solve the image restoration
problem and also compare the efficiency of Algorithm \ref{a} with the
Algorithm \ref{63} FBS, AVFB, FBMMMA and NAGA. All algorithms were written
in Matlab 2020b and an Asus Intel Core i7 laptop.

Recall that the image restoration problem can be formulated as the following
linear inverse problem:
\begin{equation}
\upsilon =Ax+b,  \label{26}
\end{equation}%
where $x\in
\mathbb{R}
^{n\times 1}$ is the original image, $A$ $\in
\mathbb{R}
^{m\times n}$ is a blurring matrix, $b\in $ $%
\mathbb{R}
^{m\times 1}$ is the additive noise and $v$ is the observed image. Also, it
is well-known that solving (\ref{26}) is equivalent to solving the convex
minimization problem%
\begin{equation}
x^{\ast }=\argmin_{x\in
\mathbb{R}
^{n}}\left\{ \frac{1}{2}\left\Vert Ax-v\right\Vert _{2}^{2}+\tau \left\Vert
x\right\Vert _{1}\right\} .  \label{27}
\end{equation}%
To solve the problem (\ref{27}), we use Theorem \ref{the}. We set $h\left(
x\right) =\frac{1}{2}\left\Vert Ax-v\right\Vert _{2}^{2}$ \ and $g\left(
x\right) =\tau \left\Vert x\right\Vert _{1}$. Then, it is easy to see that
the gradient of $h$ is $\nabla h\left( x\right) =A^{T}\left( Ax-v\right) $,
where $A^{T}$ is a transpose of $A.$ The Lipschitz constant $L$ of $\nabla h$
is computed by the maximum eigenvalues of the matrix $A^{T}A.$ In all
comparisons, we use the test images Lena, Pepper and Cameraman. At this
point, the degenerate image is obtained by adding motion blur and random
noise to the test images. In order to add the blur, we use the Matlab
function fspecial('motion', 15,60). We will try to obtain an image close to
the original image using Algorithm \ref{a}. The quality of the restored
image is measured by the signal to noise ratio (SNR) which is defined by%
\begin{equation*}
SNR=20\log \frac{\left\Vert x\right\Vert _{2}}{\left\Vert x-x_{n}\right\Vert
_{2}}\text{,}
\end{equation*}%
where $x$ and $x_{n}$ is the original image and the estimated image at
iteration $n$, respectively.

\begin{figure}[htbp]
\begin{center}
\includegraphics[width=15.0cm]{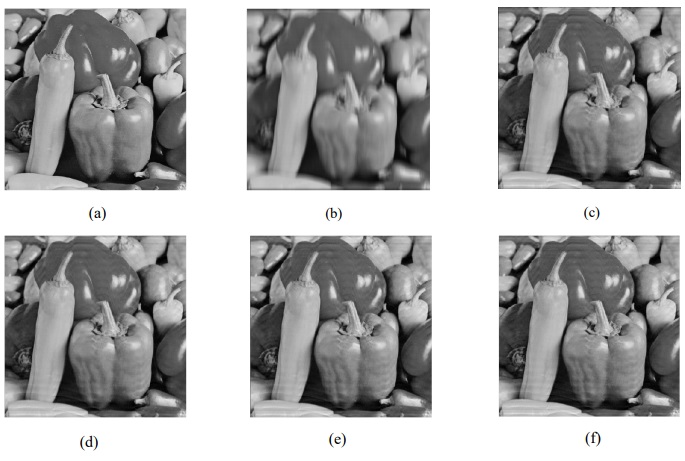}
\end{center}
\caption{(a) Pepper image (b) Blurred image (c) exclusion of viscosity and
inertial terms from Algorithm \protect\ref{a} (d) exclusion of viscosity
term from Algorithm \protect\ref{a} (e) exclusion of inertial term from
Algorithm \protect\ref{a} (f) Algorithm \protect\ref{a}}
\label{fig3b}
\end{figure}

\begin{figure}[tbph]
\begin{center}
\includegraphics[width=9.0cm]{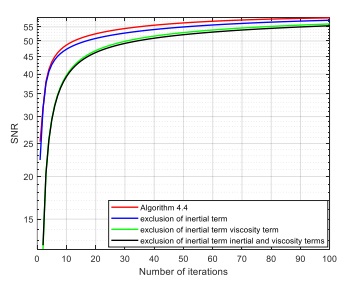}
\end{center}
\caption{ Graphic of SNR values for Pepper image}
\label{fig4b}
\end{figure}

In Theorem \ref{the}, we take $\lambda _{n}=\frac{n}{L\left( n+1\right) }$, $%
\beta _{n}=\gamma _{n}=\frac{1}{n}$, $\alpha _{n}=\frac{1}{2}$, $\eta _{n}=%
\frac{10^{20}}{n}$ and $f\left( x\right) =0.99x$. If $\theta =0.99$, then we
obtain
\begin{equation*}
\theta _{n}:=\left\{
\begin{array}{ll}
\min \left\{ 0.99,\frac{10^{20}}{n^{2}\left\Vert x_{n}-x_{n-1}\right\Vert }%
\right\} & if\text{ \ }x_{n}\neq x_{n-1}, \\
0.99, & otherwise.%
\end{array}%
\right.
\end{equation*}%
Firstly, we investigate the effectiveness of inertial and viscosity terms in
Algorithm \ref{a}. The obtained results are given in Figure \ref{fig3b},
Figure \ref{fig4b} and Table \ref{table2}.

\begin{table}[tbph]
\begin{center}
\begin{tabular}{ccccc}
\hline
&  & SNR Values &  &  \\ \hline
No.Iterations & Algorithm \ref{a}$\qquad $ & Exclusion of inertial&
Exclusion of viscosity & Exclusion of inertial and viscosity \\
$1$ & $25.292058$ & $22.368991$ & $0.000000$ & $0.000000$ \\
$5$ & $43.403404$ & $42.632692$ & $29.530821$ & $29.345449$ \\
$10$ & $48.643235$ & $47.275684$ & $39.653579$ & $39.294982$ \\
$20$ & $52.490425$ & $50.871112$ & $46.844475$ & $46.247206$ \\
$30$ & $54.397355$ & $52.757836$ & $49.890838$ & $49.194864$ \\
$40$ & $55.626238$ & $54.027441$ & $51.688664$ & $50.950629$ \\
$50$ & $56.494179$ & $54.971150$ & $52.927414$ & $52.174183$ \\
$70$ & $57.613901$ & $56.302456$ & $54.587938$ & $53.845230$ \\
$80$ & $57.977812$ & $56.788446$ & $55.183919$ & $54.458482$ \\
$90$ & $57.977812$ & $57.190757$ & $55.679869$ & $54.976793$ \\
$100$ & $58.458689$ & $57.526028$ & $56.098940$ & $55.421989$ \\ \hline
\end{tabular}%
\end{center}
\caption{SNR values for Pepper image }
\label{table2}
\end{table}

The experimental results show that the viscosity and inertial terms enable
us to obtain more quality results in image restoration.

\begin{figure}[htbp]
\begin{center}
\includegraphics[width=15.0cm]{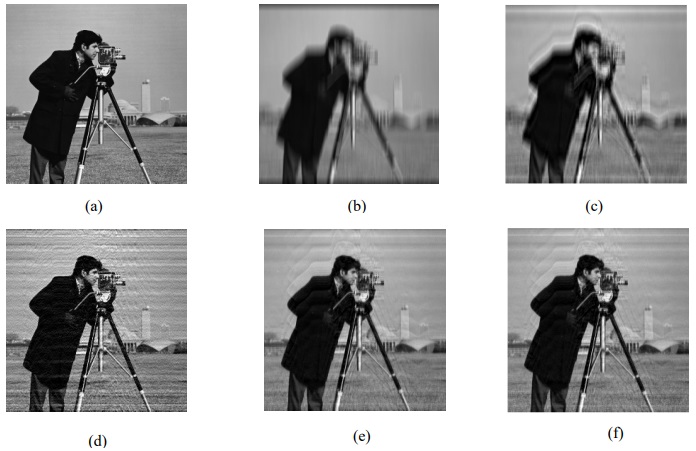}
\end{center}
\caption{(a) Cameraman image (b) Blurred image (c) NAGA (d)FBMMMA (e) AVFBA
(f) Algorithm \protect\ref{a}}
\label{fig6b}
\end{figure}

\begin{figure}[tbph]
\begin{center}
\includegraphics[width=9.0cm]{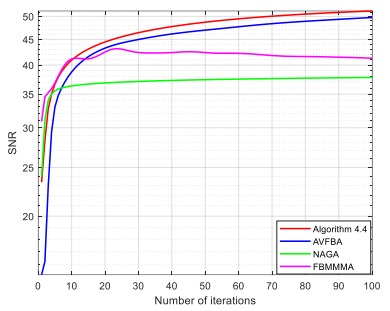}
\end{center}
\caption{ Graphic of SNR values for Cameraman image}
\label{fig7b}
\end{figure}

In what follows, we compare Algorithm \ref{a} with AVFB, FBMMMA, and NAGA
for the infinite family of nonexpansive mappings. We take the parameters as $%
\lambda _{n}=\frac{n}{L\left( n+1\right) }$, $\beta _{n}=\gamma _{n}=\rho
_{n}=\sigma _{n}=\frac{1}{n}$, $\alpha _{n}=\delta _{n}=\varphi _{n}=\frac{1%
}{2}$, $\varrho _{n}=1-\frac{1}{2}-\frac{1}{n},$ $\eta _{n}=\frac{10^{20}}{n}%
,$ $f\left( x\right) =0.7x$. In this case, the experimental results have
been given in Figure \ref{fig6b} Figure \ref{fig7b} and Table \ref{table3}.

\begin{table}[tbph]
\begin{center}
\begin{tabular}{ccccc}
\hline
&  & SNR Values &  &  \\ \hline
No.Iterations & Algorithm \ref{a} $\qquad $ & AVFBA & FBMMMA & NAGA \\
$1$ & $23.367583$ & $15.283348$ & $30.855847$ & $24.005146$ \\
$5$ & $36.649323$ & $33.006604$ & $36.84145$ & $35.540217$ \\
$10$ & $40.897499$ & $38.751068$ & $41.140872$ & $36.345638$ \\
$20$ & $44.460392$ & $43.223013$ & $42.465793$ & $36.843940$ \\
$30$ & $46.393900$ & $44.970841$ & $42.365621$ & $37.097989$ \\
$40$ & $47.709275$ & $46.126490$ & $42.353539$ & $37.269609$ \\
$50$ & $48.684256$ & $46.944487$ & $42.355574$ & $37.398947$ \\
$70$ & $50.040212$ & $48.358723$ & $41.786430$ & $37.588697$ \\
$80$ & $50.525856$ & $48.881213$ & $41.589025$ & $37.662421$ \\
$90$ & $50.924759$ & $49.328276$ & $41$.$462924$ & $37.726771$ \\
$100$ & $51.256042$ & $49.738586$ & $41.258755$ & $37.783805$ \\ \hline
\end{tabular}%
\end{center}
\caption{ SNR values for Cameraman image}
\label{table3}
\end{table}

Algorithm \ref{a} provides image restoration with higher SNR, so the
performance of image restoration of Algorithm \ref{a} is better than FBMMMA,
AVFB, and NAGA.

\begin{figure}[htbp]
\begin{center}
\includegraphics[width=15.0cm]{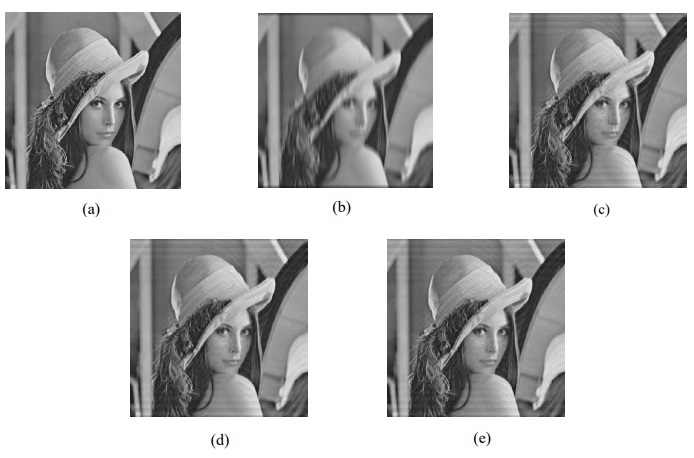}
\end{center}
\caption{ (a) Lena image (b) Blurred image (c) FBS (d) Algorithm \protect\ref%
{63} (e) Algorithm \protect\ref{a}}
\label{fig9b}
\end{figure}

\begin{figure}[tbph]
\begin{center}
\includegraphics[width=9.0cm]{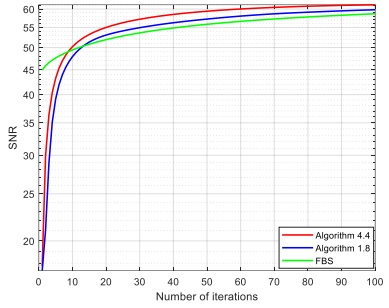}
\end{center}
\caption{ Graphic of SNR values for Lena image}
\label{fig10b}
\end{figure}

Finally, we compare Algorithm \ref{a} with Algorithm (\ref{63}) and FBS
algorithm. If we choose $\lambda _{n}=\frac{1}{L},$ $\beta _{n}=\gamma
_{n}=\delta _{n}=\psi _{n}=\phi _{n}=\frac{1}{n},\eta _{n}=\frac{10^{20}}{n}%
, $ $\theta =0.99,\alpha _{n}=\frac{1}{2},$ and $\nabla h\left( x\right)
=f\left( x\right) =0.6x$ then the experimental results are given in Figure %
\ref{fig9b}, Figure \ref{fig10b} and Table \ref{table4} .

\begin{table}[tbph]
\begin{center}
\begin{tabular}{cccc}
\hline
& SNR Values &  &  \\ \hline
No.Iteration & Algorithm \ref{a} & Algorithm \ref{63} & FBS \\
$1$ & $17.748951$ & $17.372599$ & $44.922549$ \\
$5$ & $43.22688$ & $39,.82411$ & $47.593908$ \\
$10$ & $50.193068$ & $47.884830$ & $49.431945$ \\
$20$ & $54.994313$ & $53.046121$ & $51.900830$ \\
$30$ & $57.174776$ & $54.929475$ & $53.597159$ \\
$40$ & $58.497187$ & $56.231519$ & $54.849974$ \\
$50$ & $59.389414$ & $57.219987$ & $55.820225$ \\
$70$ & $60.479787$ & $58.659111$ & $57.246811$ \\
$80$ & $60.814045$ & $59.112947$ & $57.793841$ \\
$90$ & $61.056274$ & $59.487271$ & $58.264392$ \\
$100$ & $61.228103$ & $59.805114$ & $58.673967$ \\ \hline
\end{tabular}%
\end{center}
\caption{SNR values for Lena image}
\label{table4}
\end{table}

We deduce that Algorithm \ref{a} have higher SNR than Algorithm \ref{63} and
FBS algorithm. That is, Algorithm \ref{a} has more effective in image
restoration than the other algorithms.


\begin{thebibliography}{99}
\bibitem{aoya} Aoyama K., Yasunori Y., Takahashi W., Toyoda M., On a
strongly nonexpansive sequence in a Hilbert space, J. Nonlinear Convex
Anal., \textbf{8}, (2007), 471-489.

\bibitem{in} Artsawang N., Ungchittrakool K., Inertial Mann-Type Algorithm
for a Nonexpansive Mapping to solve Monotone Inclusion and Image Restoraton
Problem, Symmetry, 12, (5) (2020) 750.

\bibitem{bau} Bauschke H.H., Combettes P. L., Convex Analysis and Monotone
Operator Theory in Hilbert Space, CMS Books in Mathematics Springer, New
York, 2011.

\bibitem{bussaban} Bussaban L., Suantai S., Kaewkhao A., A parallel
inerrtial S-iteration forward-bacward algorithm for regression and
classification problems, Carpathian. J., \textbf{36}, (2020), 35-44.

\bibitem{by} Byrne C., Aunified treatment of some iterative algorithms in
signal processing and image reconstruction, Inverse Probl., \textbf{20},
(2004), 103-120.

\bibitem{cho} Cholamjiak, P.; Shehu, Y. Inertial forward-backward splitting
method in Banach spaces with application to compressed sensing, Appl. Math.,
\textbf{64}, (2019), 409-435.

\bibitem{choo} Cholamjiak P, Suantai S. Viscosity approximation methods for
a nonexpansive semigroup in Banach spaces with gauge functions, J. Glob.
Optim.,\textbf{54,} (2012), 185-197.

\bibitem{ch} Cholamjiak W., Cholamjiak, P., Suantai, S., An inertial
forward-backward splitting method for solving inclusion problems in Hilbert
spaces, J. Fixed Point Theory Appl., (2018)
https://doi.org/10.1007/s11784-018-0526-5.

\bibitem{com} Combettes, P.L.; Wajs, V. Signal recovery by proximal
forward-backward splitting, Multiscale Model. Simul., \textbf{\ 4}, (2005),
1168-1200.

\bibitem{dong} Dong, Q., Jiang, D., Cholamjiak, P., Shehu, Y., A strong
convergence result involving an inertial forward-backward algorithm for
monotone inclusions, J. Fixed Point Theory Appl. 19 (2017) 3097--3118.

\bibitem{gobel} Goebel K., Kirk W.A, Topics in Metric Fixed Point Theory,
Cambridge Studies in Advanced Mathematics, Cambridge University Press,
Cambridge, UK, 1990.

\bibitem{gobell} Goebel K., Reich S., Uniform convexity, hyperbolic geometry
and nonexpansive mappings, Marcel Dekker, New York, NY,USA, 1984.

\bibitem{han} A. Hanjing, S. Suantai, A fast image restoration algorithm
based on a fixed point and optimization method, Mathematics 8 (2020) 378.

\bibitem{hussain} Hussain N., Ullah K, Arshad M., Fixed point approximation
of Suzuki generalized nonexpansive mappings via new faster iteration
process, J. Nonlinear Convex Anal., \textbf{19}, (2018), 1383-1393.

\bibitem{kit} Kitkuan D., Kumam P., Martifez-MoreNo J., Sitthithakerngkiet
K., Inertial viscosity forward-bakward splitting algorithm for monotone
inclusions and its application to image restoration problems, Inter. J.
Computer Math., \textbf{97}, (2019), 482-497.

\bibitem{kit1} Kitkuan D., Kumam P. and Martifez-Moreno J., Generalized
Halpern-type forward-backward splitting methods for convex minimization
problems with application to image restoration problems, Optimization,
(2019), 1-25.

\bibitem{kun} Kunrada, K.; Pholasa, N.; Cholamjiak, P. On convergence and
complexity of the modified forward-backward method involving new
linesearches for convex minimization, Math. Meth. Appl. Sci., \textbf{42},
(2019),\textbf{\ }1352-1362.

\bibitem{lion} Lions P.L., Mercier B., Splitting algorithms for the sum of
two nonlinear operators. SIAM J. Numer. Anal.,\textbf{16}, (1979), 964--979.

\bibitem{lo} Lorenz D., Pock T., An inertial forward-backward algorithm for
monotone inclusions, J. Math. Imaging Vis., \textbf{51}, (2015), 383-390.

\bibitem{main} Mainge P. E., Strong convergence of projected subgradient
methods for nonsmooth and nonstrictly convex minimization, Set Valued Anal.,
\textbf{16}, (2008), 899-912.

\bibitem{nan} Nakajo K., Shimoji K, Takahashi W., Strong convergence to a
common fixed point of families of nonexpansive mappings in banach spaces. J.
Math. Anal. Appl., \textbf{8} ,(2007), 11-34.

\bibitem{kumam} Padcharoen A., Kumam P., Fixed point optimization method for
image restoration, Thai J. of Math., \textbf{18}, (2020), 1581-1596.

\bibitem{sut} Puangpee J., Suantai S., A new accelerated viscosity iterative
method for an infinite family of nonexpansive mappings with applications to
image restoration problems, Mathematics \textbf{8}, 2020.

\bibitem{tak} Takahashi, W. Viscosity approximation methods for countable
families of nonexpansive mappings in banach spaces., Nonlinear Anal.,\textbf{%
70} ,(2009), 719-734.

\bibitem{kittt} Thong, D.V.; Cholamjiak, P. Strong convergence of a
forward-backward splitting method with a new step size for solving monotone
inclusions, Comput. Appl. Math., \textbf{38}, 2019 .

\bibitem{taka} Takahashi W., Wong N.C., Yao J.C., Two generalized strong
convergence theorems of Halpern's type in Hilbert spaces and applications,
Taiwan J. Math., \textbf{16}, (20129, 1151-1172.

\bibitem{1} Tikhonov, A.N., Arsenin V.Y., Solution of Ill-Posed Problems,
V.H. Winston: Washington, DC, USA, 1977.

\bibitem{xu} Xu, H. K., Another control condition in an iterative method for
nonexpansive mappings. Bull. Aust. Math. Soc, \textbf{65}, (2002), 109-113.

\bibitem{inertial} Verma M., Shukla K. K., A new accelerated proximal
gradient technique for regularized multittask learning framework, Pattern
Recogn. Lett. \textbf{95}, (2017), 98-103.
\end{thebibliography}
\end{document}